\documentclass[sn-mathphys-ay]{sn-jnl}% Math and Physical Sciences Author Year Reference Style
\usepackage{graphicx}%
\usepackage{multirow}%
\usepackage{amsmath,amssymb,amsfonts}%
\usepackage{amsthm}%
\usepackage{mathrsfs}%
\usepackage[title]{appendix}%
\usepackage{xcolor}%
\usepackage{textcomp}%
\usepackage{manyfoot}%
\usepackage{booktabs}%
\usepackage{algorithm}%
\usepackage{algorithmicx}%
\usepackage{algpseudocode}%
\usepackage{listings}%
\usepackage{pdflscape}
\usepackage{natbib}
\def\nbI{\ensuremath{\mathrm{1\!l}}}
\newtheorem{thm}{Theorem}%  meant for continuous numbers
\newtheorem{prop}[thm]{Proposition}% 
\newtheorem{cor}[thm]{Corollary}
\newtheorem{lem}[thm]{Lemma}

\newtheorem{rem}{Remark}%

\begin{document}

\title[Multiple change-points detection based on 
$U$-Statistics under weak dependence]{Multiple change-points detection based on 
$U$-Statistics under weak dependence }


\author[1]{\fnm{Joseph} \sur{ Ngatchou-Wandji}}\email{joseph.ngatchou-wandji@univ-lorraine.fr}

\author[2]{\fnm{Echarif} \sur{Elharfaoui}}\email{michel.harel@unilim.fr}

\author[3,4,5]{\fnm{Michel} \sur{Harel}}\email{elharfaoui.e@ucd.ac.ma}

\affil[1]{\orgdiv{ Universit\'e de Rennes (EHESP)} \orgname{\& Institut \'Elie Cartan de Lorraine}, \orgaddress{\city{Nancy}, \country{France}}}

\affil[2]{\orgdiv{Department of Mathematics, Faculty of Sciences}, \orgname{Chouaib Doukkali University}, \orgaddress{ \city{El Jadida}, \country{Morocco}}}

\affil[3]{\orgdiv{INSP\'E}, \orgname{Universit\'e de Limoges}, \orgaddress{ \city{Limoges}, \country{France}}}

 \affil[4]{\orgdiv{IMT, UMR 5219}, \orgname{Universit\'e Paul Sabatier}, \orgaddress{ \city{Toulouse}, \country{France}}}

 \affil[5]{\orgdiv{Vie-Sant\'e, UR 24134, Facult\'e de M\'edecine}, \orgname{Universit\'e de Limoges}, \orgaddress{ \city{Limoges}, \country{France}}}

\abstract{We study multiple change-points detection using
multi-samples tests based on $U$-statistics for absolutely regular observations.
Our results extend those of \cite{nga} concerned with the study of one single change-point. 
The asymptotic distributions of the test statistics under the null hypothesis and under
a sequence of local alternatives are given explicitly, and the tests are shown to be consistent. A small set of simulations is done for evaluating the performance of the tests in detecting multiple changes in the mean, variance and autocorrelation of some simple times series models.}
 


  \maketitle

\section{Introduction}

Testing for changes in the structure of data is known in statistics as a change-point problem. This is 
of great importance in a wide range of disciplines, such as industrial quality control, financial market, 
epidemic, medical diagnostics, hydrology and climatology, to name a few. 
%\smallskip

\begin{figure}
%\vspace{6pc}
   \begin{center}
      \includegraphics[width = 140mm,height=140mm]{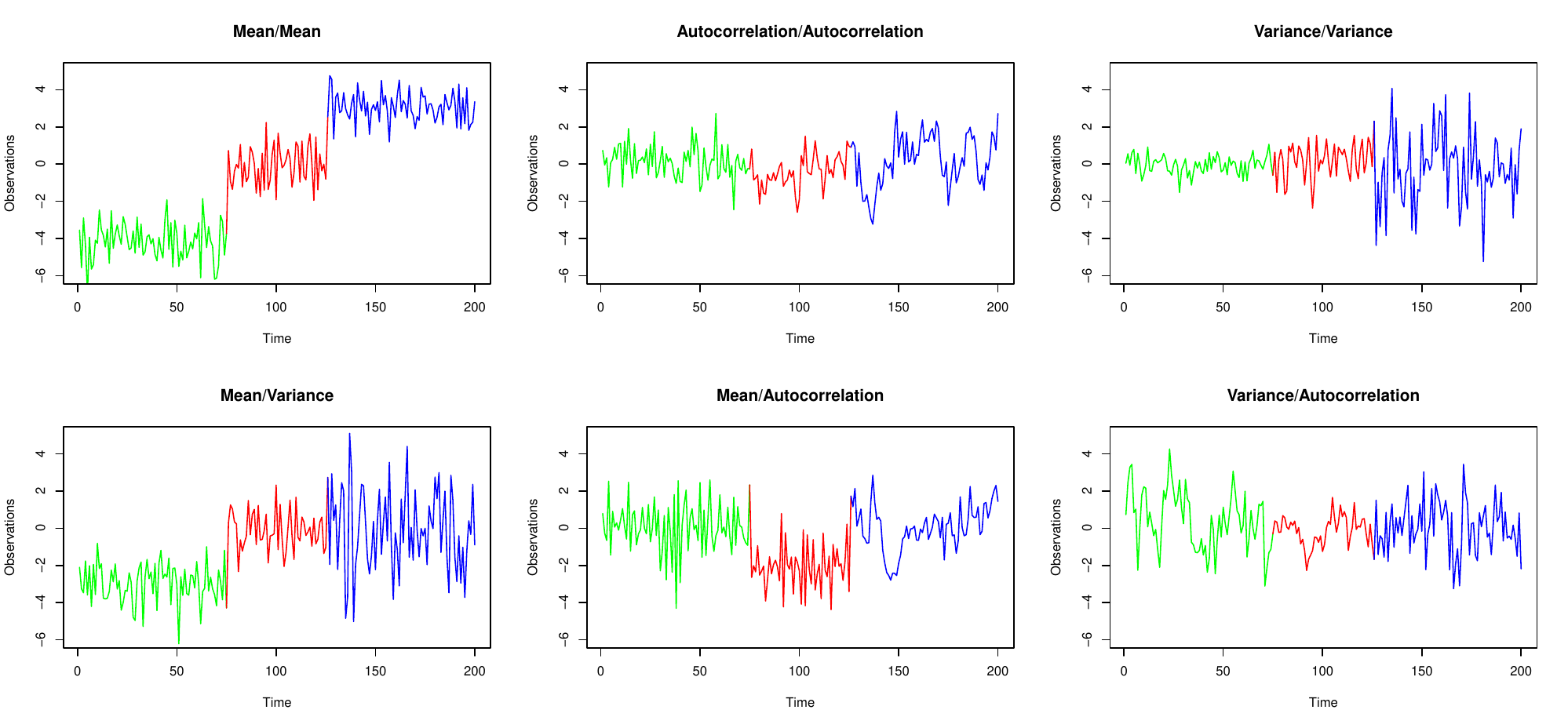}
  \end{center}
\caption{Chronograms of series with two changes. On the first row, 
the first graphic is that of a series with two changes in the mean. The second is that of a series with two changes in the autocorrolation, the last is that of a series with two changes in the variance. On the second row, the first chronogram is that of a series with one change in the mean and one in the autocorrelation. The next is that of a series with one  changes in the mean and one in the variance, and the last is that of a series with one change in the variance and another one in the autocorrelation. }
\label{f}
\end{figure}

Figure \ref{f} presents the chronograms of series with two changes, simulated from the model (\ref{mod}) in Section 4. Both changes occured either in the mean, or in the variance, or in the autocorrelation, or in the mean and the autocorrelation, or in the mean and the variance, or in the variance and the autocorrelation.
%The first chronogram is that of a series with two changes in the mean. The second is that of a series with two changes in the autocorrelation, the third is that of a series with two changes in the variance, the fourth is the chronogram of a series with one change in the mean and one in the autocorrelation, the fifth is that of a series with one  change in the mean and one in the variance, and the last is that of a series with one change in the variance and another one in the autocorrelation. 
Given such series, the change-point theory provides methods for their detection and for estimating their locations.
%\smallskip

The literature on change-point is vast. Parametric as well as nonparametric methods are used for both independent and dependent observations.  The major part of the earliest works was done in the independent context. A small sample of these are 
\cite{page1} who first, studied change-point problem, \cite{cher} who propose
test statistics for shifts detection in the mean of a normal distribution function.
This work was extended to exponential family by \cite{kan}, \cite{gard} and \cite{mac}. \cite{mat} study maximal score statistics to test for constant hazard against a change-point alternative.
\cite{hac} propose a likelihood ratio test for a change-point in a sequence
of independent exponentially distributed random variables, and prove its optimality in the sense of Bahadur.  
%\smallskip

Most of the recent papers on change-points are in the dependence context. Among them are 
\cite{vogel1} who constructs Wald-type tests for breaks detection in the trend of a dynamic time series. \cite{vogel2} who studies tests for change detection in the mean of various time series models. 
 \cite{berk}, \cite{gom0}, \cite{gom}, \cite{bardet_kw},  \cite{kengne},  \cite{bardet_k} study change detection in the parameters of various sub-classes of nonlinear heteroscedastic time series models. More recently,
\cite{ma}, \cite{mohr} and \cite{nga1, nga2} have also 
studied change detection within sub-classes of nonlinear heteroscedastic time series models.
%\smallskip

However, a substantial part of the literature is devoted to the case of one single change-point. In this paper, we stress on multiple change-points detection for weakly dependent observations containing a given number of changes. This important problem has attracted more and more attention recently. Our approach, based on $U$-statistics, is an extension of \cite{nga} who study the case of one single change.

There are various approaches in the literature dealing with $U$-statistics in change point analysis.
\cite{doring1} proposes an estimation of change-points by using $U$-statistics, while \cite{doring2} studies the convergence of such estimators. 
%Some papers dealing with change-point problem based on $U$-statistics are, among others, 
\cite{fer} estimates change location in the independence context. 
\cite{laj} compute critical values for various tests based on $U$-statistics to detect a possible change. \cite{kir} propose a general framework of sequential change-point testing procedures.
\cite{deh3} study the large-sample behavior of change-point tests in the case of short-range dependent data.
%\cite{yang}, \cite{mohr} and Yang et al. (2020). 
Additional change-point studies based on $U$-statistics are, among others \cite{aly}, \cite{gom01, gom02},  \cite{ora}, \cite{laj}, \cite{deh13} and \cite{deh17}, \cite{rack}.
% and \cite{nga}.
%\smallskip

The local power of existing change-points tests are generally not studied, in particular those of multichange-points tests. For the tests studied in this paper, not only their null distributions and consistency are studied, but also their local powers. These tests are derived from a basic process of the following form
\begin{equation*}
\mathcal Z_{n}^*(t_1,t_2,\ldots,t_k)=n^{-3/2}\sum_{l=1}^{k}\sum_{i=[nt_{l-1}]+1}^{[nt_{l}]}
\sum_{j=[nt_{l}]+1}^{[nt_{l+1}]}h(X_{i},X_{j}), 
\end{equation*}
for $\ 0\leq t_{l-1}< t_l\leq1, \ 1\leq l\leq k+1$, where $k$ stands for the number of changes in the series $X_1, \ldots, X_n$, $n$ is the sample length, $[nt_0]=1$, $[nt_{k+1}]=n$ (by convention) and $h:\mathbb{R}^{2}\rightarrow \mathbb{R}$ is a kernel function. This process is a generalization of that  defined in \cite{nga} from the one single change-point to the multiple change-points.
The study of the asymptotic properties of such processes is simplified by the use of $U$-statistics and the Hoeffding decomposition (see, eg, \cite{hoef}).
%\smallskip

The paper is organized as follows. In Section 2, we give some useful definitions and  list the main assumptions. In Section 3 we study
the asymptotic properties of our test statistics under the null hypothesis, under a
sequence of local alternatives and under fixed alternatives. Practical considerations
are presented and discussed in Section 4. Section 5 concludes our work, while the last section contains the proofs of the main results.

%\begin{center}
%\resizebox{14cm}{10cm}{\includegraphics{exemple}}
%\end{center}

\section{Definitions and assumptions}
In this section, we define our test statistics and give the main assumptions needed for the study of their behavior under the hypotheses considered. 

Let $Q$ and $R$ be cumulative distribution functions. Let $h$ be a kernel defined on $\mathbb R^2$. Define by $\theta_h(Q,R)$ the following
real number
$$\theta_h(Q,R) =\int\int h(x, y)dQ(x)dR(y).$$
For any $i = 1, 2, \ldots, n$, denote by $F_i$ be the cumulative distribution function of $X_i$. Our aim is to check the equality of the $F_i$'s. For doing this, we construct tests for testing the hypothesis $\mathcal{H}_{0}$ against the fixed alternative $\mathcal{H}_1^k$, or against the sequence of local alternatives $\mathcal{H}_{1,n}^k$
defined respectively by
\medskip

\begin{itemize}
\item[$\mathcal{H}_{0}$: ]  The distribution function of $ X_{i}
\mbox{ is } F_{1} \mbox{ for all } i\in \{1,2,\ldots,n\}$
\item[$\mathcal{H}_1^k$: ]  $\mbox{ There exist } k+2 \mbox{ real numbers }:
0=t_{0,0}<t_{0,1}<t_{0,2}<\ldots<t_{0,k}<t_{0,k+1}=1$
such that $\mbox{the cumulative distribution function of } X_{i} \mbox{ is } \dot F_l \mbox{ if } [nt_{0,l-1}] \leq i<[nt_{0,l}], \ 1\leq l\leq k+1$ and the cumulative distribution function of $ X_{n} \mbox{ is } \dot F_{k+1}$, there exist some integers $l$ and $l'$, $1\leq l<l'\leq k+1$  such that $\dot F_l \ne \dot F_{l'}$ with $\theta_h(\dot F_l, \dot F_l)\neq\theta_h(\dot F_l, \dot F_{l'})$.
%\end{itemize}
%
%
%\noindent
%We also consider the local alternatives $\mathcal{H}_{1,n}^k$ of the form
%\medskip
%
%\begin{itemize}
\item[$\mathcal{H}_{1,n}^k$: ]  $ \mbox{ There exist } k+2 \mbox{ real values } : 0=t_{0,0}<t_{0,1}<t_{0,2}<\ldots<t_{0,k}<t_{0,k+1}=1$ such that $F_{[nt_l]+1}=F_{[nt_l]+2}=\ldots=F_{[nt_{l+1}]}=F_{l+1}^{(n)}$, $0 \leq l \leq k$, the distribution function of $X_n$ is $F_{k+1}^{(n)}$, there exists some $l$, $1\leq l\leq k+1$  such that $\theta_h(F_1,F^{(n)}_l)= \theta_h(F_1,F_1)+n^{-1/2}[A_l+o(1)]$, where $A_l$ is some constant.
\end{itemize}
\smallskip

\begin{rem} As example of local alternatives $\mathcal{H}_{1,n}^k$, we can take those defined by $F^{(n)}_l(x) = F_1(x+n^{-1/2}\gamma_l)$, for some $\gamma_l \in \mathbb R$, together with a 
twice differentiable kernel function $h$ satisfying $\int\int(\partial h(x, y)/\partial y)dF_1(x)dF_1(y) < \infty$ and $|\partial^2h(x, y)/\partial^2y| < C$ for some $C>0$. These can be checked easily by an 
application of the Taylor-Young formula. The computations with this example, yield, for any $l$,
$$A_l = -\gamma_l\int\int\frac{\partial h}{\partial u}(x, u)dF_1(x)dF_1(u).$$
\end{rem}
\bigskip

For solving our testing problem we can construct Kolmogorov-Smirnov (KS) or Cram\'er-von Mises (CV) type-tests based respectively on the following KS and CM statistics 
\begin{equation} \label{eq1}
\mathcal T_{1,n} = \max_{1< m_1<m_2<\ldots<m_k<n}\left|n^{-3/2}\sum^k_{l=1}\sum_{i=m_{l-1}}^{m_l}\sum_{j=m_l+1}^{m_{l+1}}\left [h(X_i,X_j)
-\theta_h(\widehat{F}_1,\widehat{F}_1)\right]\right|
\end{equation}

\begin{equation} \label{eq2}
\mathcal T_{2,n} =\frac{1}{n^k}\sum_{1< m_1<m_2<\ldots<m_k<n}\left\{n^{-3/2}\sum^k_{l=1}\sum_{i=m_{l-1}}^{m_l}\sum_{j=m_l+1}^{m_{l+1}}\left[h(X_i,X_j)
-\theta_h(\widehat{F}_1,\widehat{F}_1)\right]\right\}^2
\end{equation}
where by convention $m_0=1$ and $m_{k+1}=n$ and $\widehat{F}_1$ stands for any consistent estimator
of $F_1$, a simple example being the associated empirical cumulative distribution function. 
\medskip

\noindent
Denote by $[x]$ the integer part of any real number $x$ and by $\Theta_k$ the subset of $[0,1]^k$ defined by 
$$\Theta_k=\{t=(t_1,\ldots,t_k) \in [0,1]^k; 0<t_1< \ldots< t_k<1\}.$$
Noting that for any $k\in\{1,\ldots,n-1\}$,
there exists $t_l^* \in [0, 1] \ (1\leq l\leq k)$ such that $k = [nt_l^*]$, one can write, asymptotically,
\begin{eqnarray*}
\mathcal T_{1,n} &=& \sup_{t_1,\ldots,t_k \in [0,1]}| \mathcal Z_{n}(t_1,\ldots,t_k)| \\
\mathcal T_{2,n} &=& \int_0^1 \ldots \int_0^1 \mathcal Z_{n}^2 (t_1,\ldots,t_k)dt_1 \ldots dt_k, 
\end{eqnarray*}
where $(\mathcal Z_{n}(\cdot))$ stands for the sequence of stochastic processes defined for all  
 $t=(t_1,\ldots,t_k) \in [0,1]^k$, by 
\begin{equation} \label{eq3}
\mathcal Z_{n}(t_1,\ldots,t_k)=
\left \{
\begin{array}{lc} \displaystyle 
n^{-3/2}\sum_{l=1}^{k}\sum_{i=[nt_{l-1}]+1} ^{[nt_{l}]}
\sum_{j=[nt_{l}]+1}^{[nt_{l+1}]}\left\{h(X_{i},X_{j})-\theta_h(F_1,F_1)\right\},& t \in \Theta_k\\
0& t \notin \Theta_k, 
\end{array}
\right.
\end{equation}
with the convention that $[nt_0]=1$ and $[nt_{k+1}]=n$.
\bigskip

For any
%$0\leq a < b <c \leq 1$, define the following $U$-statistic $U_{n,a,b,c}$ with 
kernel $h$,  noting that $\dot F_1=F_1$, define the following functions used later in the Hoeffding decomposition :
%
%\begin{equation*}
%U_{n,a,b,c} =\frac{1}{([nc]-[nb])([nb]-[na])} \sum_{i=[na]+1}^{[nb]} \sum_{j=[nb]+1}^{[nc]}h(X_i,X_j)
%\end{equation*}
%
\begin{equation*}
h_{F_1,1}(x) =\int h(x, y)dF_1(y)-\theta_h(F_1,F_1), \ x \in \mathbb R
\end{equation*}
\begin{equation*}
h_{F_1,2}(y) = \int h(x, y)dF_1(x)-\theta_h(F_1, F_1), \ y \in \mathbb R
\end{equation*}
\begin{equation*}
h_{\dot F_l,1}(x) = \int h(x, y)d\dot F_l(y)-\theta_h(\dot F_{l-1},\dot F_l), \ l=2, \ldots, k, \ x \in \mathbb R
\end{equation*}
\begin{equation*}
h_{\dot F_l,2}(y) = \int h(x,y)d\dot F_l(x)-\theta_h(\dot F_{l-1},\dot F_l), \ l=2, \ldots, k, \ y \in \mathbb R
\end{equation*}
\begin{equation*}
g_{F_1,F_1}(x,y) = h(x, y)-h_{F_1,1}(x)-h_{F_1,2}(y) + \theta_h(F_1,F_1), \ x,y \in \mathbb R
\end{equation*}

\begin{equation*}
g_{\dot F_{l-1},\dot F_l}(x,y) = h(x,y)-h_{\dot F_l,1}(x)-h_{\dot F_l,2}(y) + \theta_h(\dot F_{l-1},\dot F_l), \ l=2, \ldots, k, \ x,y \in \mathbb R.
\end{equation*}
%Consider the Hoeffding decomposition of $U_{n,a,b,c}$ under $\mathcal{H}_{0}$ : 
%
%\begin{equation} \label{eq4}
%U_{n,a,b,c} = \theta_h(F_1, F_1) + U^{(1)}_{n,F_1,1} + U^{(1)}_{n,F_1,2} + U^{(2)}_{n,a,b,c},
%\end{equation}
%where
%
%$$U^{(1)}_{n,F_1,1} = \frac{1}{([nc]-[nb])}\sum_{i=[na]+1}^{[nb]} h_{F_1,1}(X_i)$$
%
%$$U^{(1)}_{n,F_1,2} = \frac{1}{([nb]-[na])}\sum_{i=[nb]+1}^{[nc]} h_{F_1,2}(X_i)$$
%
%\begin{eqnarray*}
%U^{(2)}_{n,a,b,c} & = & \frac{1}{([nc]-[nb])([nb]-[na])}
%\sum_{i=[na]+1}^{[nb]}\sum_{j=[nb]+1}^{[nc]} g_{F_1,F_1} (X_i,X_j).
%\end{eqnarray*}
%\smallskip

Also, we define the following real numbers
\begin{equation*}
\sigma_{p,r}=\mathbb{E}\big[h_{F_1,p}(X_1)h_{F_1,r}(X_1)\big]
+2\sum_{i=1}^{\infty}\mathbb{C}\mbox{ov}\big(h_{F_1,p}(X_1), h_{F_1,r}(X_{i+1})\big), \ p,r=1,2.
\end{equation*}
%\medskip

Now, we recall from \cite{har2} that a non-necessarily stationary triangular sequence $\{\mathcal{V}_{ni}, \ 1 \leq i \leq n, \ n \geq 1\}$ is absolutely regular if, as $m\rightarrow\infty$,
\begin{equation*}
\beta(m)=\sup_{n\in\mathbb{N}}\sup_{m\leq n}\max_{1\leq j\leq n-m}
\mathbb{E}\left\{\sup_{A\in\mathcal{A}_{n,j+m}^\infty}\left|\mathbb{P}(A\mid\mathcal{A}_{n,0}^j)-\mathbb{P}(A)\right| \right\} \longrightarrow0,
\end{equation*}
with $\mathcal{A}^j_{n,i}$ standing for the $\sigma$-algebra generated by $\mathcal{V}_{ni}, \ldots, \mathcal{V}_{nj}, \ i<j, \ i, j \in \mathbb{N}\cup\{\infty\}$.
\bigskip

%\noindent
In the remaining, we make the following main assumptions :\smallskip

\begin{itemize}
\item[(A1) ] The sequence $(X_i)_{i\in\mathbb{N}}$ is stationary and absolutely regular with rate
\begin{equation} \label{bet}
\beta(n)=\mathcal{O}(\tau^n), \ \ 0<\tau<1.
\end{equation}
%(A2) $\{X_i\}_{i\in\mathbb{N}}$ .\\
\item[(A2) ] For any $l=1, \ldots, k$,  consider $(Y^l_i)_{1\leq i\leq n}$ a sequence of stationary and absolute regular random variables
with rate (\ref{bet}) and cumulative distribution function $\dot F_l$.

For any $i,j\in \mathbb{N}$, the absolute
regular dependence between $Y_{i}^{l}$ and $Y_{j}^{l}$ is the same
as the dependence between $X_{i}$ and $X_{j}$. The random variables
$Y^1_i$'s have the same law as the random variables $X_i$.
\item[(A3) ] For any $l=1, \ldots, k$ consider $(Y^l_{ni})_{1\leq i\leq n,n\geq1}$ an array of stationary and absolute regular random variables with cumulative distribution function $F_l^{(n)}(x) = F_1(x + n^{-1/2}\gamma_l)$.

We assume the cumulative distribution functions $F^{(n)}_{l,ij}$ and $F^{*(n)}_{l,ij}$ of the $(Y^l_{ni}, Y^l_{nj})$'s and $(X_i, Y^l_{nj})$'s respectively satisfy
\begin{equation*}
\lim_{n\rightarrow\infty}F^{(n)}_{l,ij}(x, y) = F_{ij}(x, y) \mbox{ and }
\lim_{n\rightarrow\infty}F^{*(n)}_{l,ij}(x, y)= F_{ij}(x, y), \ 1 \leq i < j \leq n,
\end{equation*}
where $F_{ij}$ is the cumulative distribution function of $(X_i,X_j)$.
\end{itemize}

\begin{rem}
%The geometrical mixing rate is assumed by convenience. The current 
%results can be established as well for arithmetical mixing rates to be found.
Assumption (A1) is meet by  ARMA, GARCH and some other usual nonlinear time series as indicated for instance in Example 2.2 of \cite{sch}. Assumptions (A1) and  (A2) allow for the study of the behavior of the test statistic under the alternatives, while assumptions  (A1)  and (A3) allow for the study under the local alternatives defined by the sequence $(\mathcal H_{1,n}^k)$. In particular (A3) implies that $\mathcal H_0$ and  $(\mathcal H_{1,n}^k)$ are contiguous (see, e.g. \cite{lecam}).
\end{rem}
\section{Asymptotics}
This section is devoted to the statement of our theoretical results. The main proofs are postponed to the last section.
\subsection{Results for general kernels}
In this subsection, we state the results for general kernels $h$. These results show that, in general, the test statistics considered are not asymptotic distribution free. However, we explain how to approximate the limiting  distribution of  the Cram\'er-type test statistic presented later. \bigskip

\begin{thm} \label{th1}
Assume that (A1)-(A2) hold. Then, under $\mathcal{H}_{0}$, if for some
$\delta >0$
$$\max\left\{\sup_{ij}\mathbb{E}\left(\left\vert
h(X_{i},X_{j})\right\vert^{2+\delta}\right), \int\int_{\mathbb{R}^{2}}\left\vert
h(x,y)\right\vert^{2+\delta}dF_{1}(x)dF_{1}(y)\right\}<\infty,$$
%
%and the absolute regularity condition (\ref{bet}) is satisfied, 
then for any $p,r= 1, 2$, $\sigma_{pr}<\infty$. 
\medskip

If in addition $\sigma_{pr}>0$,
$1\leq p,r\leq2$, then $\{\mathcal Z_{n}(t_1,t_2,\ldots,t_k);
t_1,t_2, \ldots, t_k \in [0,1]\}_{n\in \mathbb{N}}$ converges
in distribution towards the process $\mathcal Z(\cdot)$ defined for any $t=(t_1, \ldots, t_k) \in [0,1]^k$ by 
\begin{eqnarray*}
&&\mathcal Z(t_1,t_2,\ldots,t_k)=\\
&&\left\{
\begin{array}{lc} \displaystyle 
\sum_{l=1}^{k}\left\{
(t_{l+1}-t_{l})\left[W_1(t_{l})-W_1(t_{l-1}) \right]+(t_{l}-t_{l-1})\left[W_2(t_{l+1})-W_2(t_{l})\right]\right\},& t\in \Theta_k, \\
0 & t \notin \Theta_k, 
\end{array} \right.
\end{eqnarray*}
%
%and $\mathcal Z(t_1,t_2,\ldots,t_k)=0$ otherwise. 
where by convention $t_0=0$, $t_{k+1}=1$, $\{(W_1(t),W_2(t))\}_{0\le t \le 1}$, is a two-dimensional zero-mean Brownian motion with covariance kernel matrix with entries 
%\break
$\mathbb{C}ov\big(W_p(s),W_r(t)\big)=\min(s,t)\sigma_{pr}$, $p,r= 1, 2$, $s,t \in [0,1]$.
\end{thm}
\bigskip

\noindent
{\bf Proof} -
See Appendix.
\bigskip

\begin{rem}
The covariance kernel $\Gamma$ of the limiting process $\mathcal Z (\cdot)$ is given for
any $s=(s_1, \ldots,s_k), \ t=(t_1, \ldots,t_k) \in \mathbb R^k$ 
  by $\Gamma(s,t)=0$ if $s=(s_1, \ldots,s_k), \ t=(t_1, \ldots,t_k) \notin \Theta_k$ and for 
$s=(s_1, \ldots,s_k), \ t=(t_1, \ldots,t_k) \in \Theta_k$, by 
\begin{eqnarray}
&&\Gamma(s,t)=\nonumber\\
&&\sigma_{11}\sum_{l=1}^k  (t_{l+1}-t_l)(s_{l+1}-s_l)(t_{l}\wedge s_{l}-t_{l}\wedge s_{l-1}-t_{l-1}\wedge s_{l} + t_{l-1}\wedge s_{l-1}) \nonumber\\
&+& 
\sigma_{12}\sum_{l=1}^k (t_{l+1}-t_l)(s_{l}-s_{l-1})
(t_{l}\wedge s_{l+1}-t_{l}\wedge s_{l}-t_{l-1}\wedge s_{l+1}+t_{l-1}\wedge s_{l}) \nonumber\\
&+& 
\sigma_{12}\sum_{l=1}^k (t_{l}-t_{l-1}) (s_{l+1}-s_l)
(t_{l+1}\wedge s_{l}-t_{l+1}\wedge s_{l-1}-t_{l}\wedge s_{l}+t_{l}\wedge s_{l-1}) \nonumber\\
&+& 
\sigma_{22}\sum_{l=1}^k (t_{l}-t_{l-1}) (s_l-s_{l-1})(t_{l+1}\wedge s_{l+1}-t_{l+1}\wedge s_{l}-t_{l}\wedge s_{l+1} + t_{l}\wedge s_{l}) \nonumber\\
&+& 
2\sigma_{11}\sum_{l=1}^k \sum_{j \le l} (t_{l+1}-t_l)(s_{j+1}-s_j)
(t_{l}\wedge s_{j}-t_{l}\wedge s_{j-1}-t_{l-1}\wedge s_{j} + t_{l-1}\wedge s_{j-1})
\nonumber\\
&+& 
2\sigma_{12}\sum_{l=1}^k \sum_{j \le l}  (t_{l+1}-t_l)(s_{j}-s_{j-1})
(t_{l}\wedge s_{j+1}-t_{l}\wedge s_{j}-t_{l-1}\wedge s_{j+1}+t_{l-1}\wedge s_{j})
\nonumber\\
&+& 
2\sigma_{12}\sum_{l=1}^k \sum_{j \le l} (t_{l}-t_{l-1})(s_{j+1}-s_{j})
(t_{l+1}\wedge s_{j}-t_{l+1}\wedge s_{j-1}-t_{l}\wedge s_{j}+t_{l}\wedge s_{j-1})
\nonumber\\
&+& 
2\sigma_{22}\sum_{l=1}^k \sum_{j \le l}  (t_{l}-t_{l-1})(s_{j}-s_{j-1})
%\nonumber
%\\&& \times
(t_{l+1}\wedge s_{j+1}-t_{l+1}\wedge s_{j}-t_{l}\wedge s_{j+1} + t_{l}\wedge s_{j}).
\end{eqnarray}
\end{rem}
\begin{thm} \label{th2}
Assume (A1) and (A3) hold, $h$ is twice differentiable with bounded
second-order derivatives $\partial^2h(x, y)/\partial x\partial y$, and the integral 
$\int\int(\partial h(x, y)/\partial y)dF_1(x)dF_1(y)$ is finite. Then, under $\mathcal{H}_{1,n}^k$, if there exists $\delta > 0$ such that for any $1\leq l\leq k+1$,
\begin{equation*}
\sup_{1\leq i,j\leq n}\mathbb{E}\left(\left|h(X_i,X_j)\right|^{2+\delta}\right), \ \
\sup_{n\geq1}\sup_{i,j}\mathbb{E}\left(\left|h(Y_{ni}^l, Y_{nj}^{l})\right|^{2+\delta}\right),
\end{equation*}
\begin{equation*}
\sup_{n\geq1}\sup_{1\leq i,j\leq n}\mathbb{E}\left(\left|h(X_i, Y_{nj}^{l})\right|^{2+\delta}\right), \ \
\int\int_{\mathbb{R}^2}\left|h(x, y)\right|^{2+\delta}dF_1(x)dF_1(y),
\end{equation*}
\begin{equation*}
\sup_{n\geq1}\int\int_{\mathbb{R}^2}\left|h(x, y)\right|^{2+\delta}dF^{(n)}_l(x)dF^{(n)}_{l}(y), \ \
\sup_{n\geq1}\int\int_{\mathbb{R}^2}\left|h(x, y)\right|^{2+\delta}dF_1(x)dF^{(n)}_{l}(y)
\end{equation*}
are finite,  if for any $p,r= 1, 2$, $\sigma_{pr}>0$,
then the sequence of processes $\{\mathcal Z_{n}(t_1,t_2,\ldots,t_k);$
$t_1, t_2,\ldots, t_k \in [0,1]\}_{n\in \mathbb{N}}$
converges in distribution towards a Gaussian process $\widetilde{\mathcal Z}(\cdot)$ 
%with mean $\sum_{l=1}^k(t_{l+1} -t_l)(t_l-t_{l-1})A_l$ 
with representation given for any $t=(t_1,\ldots,t_k) \in [0,1]^k$ by
\begin{equation*}
\widetilde{\mathcal Z}(t_1,t_2,\ldots,t_k)= 
\left \{
\begin{array}{lc}
\mathcal Z(t_1,t_2,\ldots,t_k)+ \displaystyle \sum_{l=1}^k(t_{l+1}-t_l)(t_l-t_{l-1})A_l, &
t \in \Theta_k\\
0 & t \notin \Theta_k,
\end{array}
\right.
\end{equation*}
where by convention $t_0=0$, $t_{k+1}=1$ and $\{\mathcal Z(t_1,t_2,\ldots,t_k), \ {t_1, t_2,\ldots, t_k \in [0,1]}\}$ is the zero-mean Gaussian process defined in Theorem \ref{th1} and the constants $A_l$ are invoked in $\mathcal H_{1,n}^k$.
\end{thm}
\bigskip

\noindent 
{\bf Proof} -
See Appendix.
\bigskip

\begin{cor} \label{cor1}
Assume that the cumulative distribution function of $X_i$ is $\dot F_l$ for 
$[nt_{0,l-1}] \leq i < [nt_{0,l}]$, $1 \leq l \leq k +1$ and that of $X_n$ is $\dot F_{k+1}$.
% for all $x \in\mathbb{R}$. 
Assume 
there exist constants $B_l$, $1\leq l\leq k+1$ such that
$F_l^{(n)}(x)=F_1(x+n^{-1/2}B_l)$ for all $x \in \mathbb{R}$, and that the kernel function $h$ is twice differentiable with 
$\int\int(\partial h(x,y)/\partial y)dF_1(x)d\dot F_l(y)<\infty$ and
$\partial^2h(x,y)/ \partial^2y$ bounded.
%
% then we are under the alternative $\mathcal{H}_{1,n}^k$.
Then Theorem 2 holds.
\end{cor}
\bigskip

\noindent 
{\bf Proof} -
It suffices to check the assumptions of Theorem \ref{th2}.
\bigskip

\begin{thm} \label{th3}
Assume (A1)-(A3) hold and that under $\mathcal{H}_1^k$, the integrability conditions
in Theorem \ref{th2} hold. Then, as $n$ tends to infinity, the sequence of processes \break 
$n^{-1/2} \mathcal Z_{n}^*(t_{1},t_{2},\ldots,t_{k})$  converges 
almost surely in $D([0,1]^k)$ to  either of the following limits 
\begin{equation} \label{z_1*}
%\left \{
\begin{array}{l}
%&&\nbI(t_1, \ldots, t_k \in [t_{0,l-1},t_{0,l}]) \times \nonumber \\
(i) \ \sum_{l=1}^{k}\left [\theta_h(\dot F_l,\dot F_l)(t_l-t_{0,l-1})(t_{0,1}-t_l)+\theta_h(\dot F_l,\dot F_{l+1})(t_l-t_{0,l-1})(t_{0,l+1}-t_{0,l})\right],
\\
 \ \ \ \ \ \mbox{if } t_{0,l-1}\le t_l \le t_{0,l}, \ l=1, \ldots, k, (t_{1},t_{2},\ldots,t_{k}) \in \Theta_k ; 
\\
%&+& \nonumber\\
%&&\nbI(t_1, \ldots, t_k \in [t_{0,l},t_{0,l+1}]) \times \nonumber \\
(ii) \  \sum_{l=1}^{k}\left[\theta_h(\dot F_{l+1},\dot F_{l+1})(t_{l}-t_{0,l-1})(t_{0,l+1}-t_{l})+\theta_h(\dot F_l,\dot F_{l+1})t_{0,l}(t_{0,l+1}-t_{l})\right], \\
%\nbI(t_l \in [t_{0,l},t_{0,l+1}])
 \ \ \ \ \  \mbox{if } t_{0,l}\leq t_{l}\leq t_{0,l+1}, \ l=1, \ldots, k, (t_{1},t_{2},\ldots,t_{k}) \in \Theta_k ;\\
(iii) \ 0, \mbox{  if } (t_{1},t_{2},\ldots,t_{k}) \notin \Theta_k
\end{array} 
%\right.
\end{equation}
 where by convention $t_{0,0}=0$.
\end{thm}
\bigskip

\noindent
{\bf Proof} - 
See Appendix.
\bigskip

\noindent
From Theorem \ref{th3}, we easily deduce the following corollary.\bigskip

\begin{cor} \label{cor2}
Under $\mathcal{H}_{0}$ and the conditions of Theorem \ref{th1}, as $n$ tends to infinity, almost surely, 
\begin{equation*}
n^{-1/2}\mathcal Z_{n}^*(t_{1},t_{2},\ldots,t_{k}) \longrightarrow
\theta_h(F_{1},F_{1})\sum_{l=1}^{k}t_{l}(t_{0,l+1}-t_{l}), \ (t_1, \ldots, t_k) \in \Theta_k.
\end{equation*}
\end{cor}
 Let $\lambda$ be any square integrable function on $[0,1]^k$. Define the integral operator $\mathcal K_\Gamma$ by 
\begin{equation} \label{op} \mathcal K_\Gamma[\lambda(\cdot)]= \int_{[0,1]^k} \Gamma(\cdot,s) \lambda(s)ds,
\end{equation}
where we recall that $\Gamma$ is the covariance kernel of the process $\mathcal Z$ defined in Theorem \ref{th1}. \bigskip

\begin{thm} \label{th4}
Assume that the assumptions of Theorem \ref{th2} hold. Let $(\mathcal Z(s) : s\in [0,1]^k)$ be the limiting process defined in Theorems \ref{th1} and \ref{th2}. Then 
\begin{itemize}
\item[i- ] Under $\mathcal H_0$, as $n$ tends to infinity, one has the following convergence in distribution, 
%\begin{eqnarray}
$$\mathcal T_{1,n} \longrightarrow  \sup_{s \in [0,1]^k}\left | \mathcal Z(s) \right | \ \mbox{ and } \ \mathcal T_{2,n} \longrightarrow \int_{[0,1]^k} \mathcal Z^2(s) ds = \sum_{j \ge 1} \zeta_j \chi_j^2,$$
where the $ \chi_j^{2}$'s are iid chi-square random variables with one degree of freedom and the $\zeta_j$'s are standing for the eigenvalues of the operator $\mathcal K_\Gamma$.
%
%\end{eqnarray}
\item[ii- ] Under $\mathcal H_{1,n}^k$, as $n$ tends to infinity, one has the following convergence in distribution, 
%\begin{eqnarray}
$$\mathcal T_{1,n} \longrightarrow 
\sup_{s \in [0,1]^k}\left | \widetilde{\mathcal Z}(s) \right | \ \mbox{ and } \ \mathcal T_{2,n} \longrightarrow  \int_{[0,1]^k} \widetilde{\mathcal Z}^2(s) ds = \sum_{j \ge 1} \zeta_j \chi_j^{*2},$$
where the $ \chi_j^{*2}$'s are iid non-central chi-square random variables with one degree of freedom and non-centrality parameters $\rho_j^2\zeta_j^{-1}$
and 
$$\rho_j=\sum_{l=1}^kA_l\int_0^1 \ldots \int_0^1 (s_{l+1}-s_l) (s_l -s_{l-1})
g_j(s_1, \ldots,s_k)\nbI((s_1,\ldots,s_k) \in \Theta_k)  ds_1 \ldots ds_k,$$
in which the $g_j$'s stand for the eigenfunctions of the integral operator $\mathcal K_\Gamma$ associated with the eigen-value $\zeta_j's$.
%\end{eqnarray}
\item[iii- ] Under $\mathcal H_1^k$, as $n$ tends to infinity, one has the following convergence in probability, 
%\begin{eqnarray}
$$\mathcal T_{1,n} \longrightarrow \infty, \ \ \ \mathcal T_{2,n} \longrightarrow \infty.$$
%\end{eqnarray}
\end{itemize}
\end{thm}
\bigskip

\noindent
{\bf Proof} - 
$i$- It is immediate from Theorem \ref{th1} that 
$\mathcal T_{1,n}$ and $\mathcal T_{2,n}$ converge in distribution respectively to 
$\sup_{t \in \Theta_k} |\mathcal Z(t)|$ and $\int_{\Theta_k} \mathcal Z^2(t)dt$.
\bigskip

We prove now that  the distribution of $\int_{\Theta_k} \mathcal Z^2(t)dt$ is that of  a sum of weighted iid chi-square distribution with one degree of freedom. 

From usual arguments, it is easy to see that the integral operator defined by (\ref{op}) admits eigenvalues $\zeta_1\ge \zeta_2\ge \ldots \ge 0$ with associated eigenfunctions $g_1,g_2, \ldots$ forming an orthonormal basis of $L^2(\Theta_k)$, the set of square integrable functions on $\Theta_k$. Then the zero-mean Gaussian process $\mathcal Z$ as a function in $L^2(\Theta_k)$, has the Karhunen-Lo\`eve representation 
$$\mathcal Z(t)=\sum_{j \ge 1} G_j g_j(t), t \in \Theta_k$$
with the independent random variables $G_j$'s defined as $G_j=\int_{\Theta_k} \mathcal Z(t)g_j(t)d t \sim \mathcal N(0, \zeta_j)$. It results from this that in distribution 
$$\int_{\Theta_k} \mathcal Z^2(t)dt=\sum_{j \ge 1}  \zeta_j \chi_j^2,$$
where the $\chi_j^2$'s are iid chi-square random variables with one degree of freedom.
\bigskip

\noindent
 $ii$- Here also, it is immediate from Theorem \ref{th2} that
$\mathcal T_{1,n}$ and $\mathcal T_{2,n}$ converge in distribution respectively to 
$$\sup_{(t_1, \ldots,t_k)  \in \Theta_k} |\mathcal Z(t_1,\ldots,t_k)+\sum_{l=1}^k(t_{l+1}-t_l)(t_l-t_{l-1})A_l| \ \mbox{ and }$$
$$\int_{[0,1]^k} \nbI((t_1,\ldots,t_k) \in \Theta_k) |\mathcal Z(t_1,\ldots,t_k)+\sum_{l=1}^k(t_{l+1}-t_l)(t_l-t_{l-1})A_l|^2 dt_1 \ldots dt_k.$$
\smallskip

\noindent
For the same reasons as above, one has the decomposition 
$$ \mathcal Z(t)+\sum_{l=1}^k(t_{l+1}-t_l)(t_l-t_{l-1})A_l=\sum_{j \ge 1} \widetilde G_j g_j(t), t=(t_1, \ldots, t_k) \in \Theta_k, $$
with the independent random variables $\widetilde G_j$'s defined as $\widetilde G_j=\int_{\Theta_k} \mathcal Z^2(t)dt \sim \mathcal N(\rho_j, \zeta_j)$, where  
$$\rho_j=\sum_{l=1}^kA_l\int_0^1 \ldots \int_0^1 (s_{l+1}-s_l) (s_l -s_{l-1})
g_j(s_1, \ldots,s_k)\nbI((s_1,\ldots,s_k) \in \Theta_k) ds_1 \ldots ds_k.$$

\noindent
 It follows from this that, in distribution, 
\begin{eqnarray*}
&&\int_0^1 \ldots \int_0^1 \nbI((t_1,\ldots,t_k) \in \Theta_k) | \mathcal Z(t_1, \ldots,t_k)+\sum_{l=1}^k(t_{l+1}-t_l)(t_l-t_{l-1})A_l|^2dt_1 \ldots dt_k\\ &&=\sum_{j \ge 1} \zeta_j \chi_j^{*2},
\end{eqnarray*}
where the $\chi_j^{*2}$'s are non-central iid chi-square random variables with one degree of freedom and non-centrality parameter $\rho_j^2 \zeta_j^{-1}$.
\bigskip

\noindent
 $iii$- The last part follows easily from Theorem \ref{th3}. 
\bigskip

\subsection{Results for some particular kernels}
For some particular kernels, the results of the above subsection simplify, and, can lead to the construction of asymptotic distribution free test statistics whose distributions can be approximated more easily than those of the preceding more general test statistics.
\bigskip

\noindent
Denote by $\sigma^2$ the long-run variance of the series $(h_{F_1,1}(X_i))$ : 
$$\sigma^2=\mathbb{E}\big\{[h_{F_1,1}(X_1)]^2\big\}
+2\sum_{i=1}^{\infty}\mathbb{C}\mbox{ov}\big(h_{F_1,1}(X_1), h_{F_1,1}(X_{i+1})\big).$$
%This number is a long-run variance.
\medskip

We have the following result.\bigskip

\begin{thm} \label{th5}
%$\bullet$
If the kernel $h$ is such that for all $l=1, \ldots,k$ its associated 
$h_{\dot F_l,1}$ and $h_{\dot F_l,2}$ satisfy
$h_{\dot F_l,1}(x)=-h_{\dot F_l,2}(x)$, then 
\begin{itemize}
\item[(i) -] Under the assumptions of Theorem \ref{th1},  under $H_0$, 
$\mathcal Z_n(t_1,t_2,\ldots,t_k)$ converges weaklyto $\mathcal Z_0(t_1,t_2,\ldots,t_k)$ with representation 
\begin{eqnarray*}
\mathcal Z_0(t_1,t_2,\ldots,t_k)&=&\left\{
\begin{array}{lc}
\sigma \mathcal B(t_1,t_2,\ldots,t_k),&  (t_1,t_2,\ldots,t_k)\in \Theta_k\\
0, & (t_1,t_2,\ldots,t_k) \notin \Theta_k, 
\end{array} \right.
\end{eqnarray*}
 where for all $(t_1,t_2,\ldots,t_k) \in \Theta_k$, 
$$\mathcal B(t_1,t_2,\ldots,t_k)=\sum_{l=1}^{k}\left\{
(t_{l+1}-t_{l})\left[W_0(t_{l})-W_0(t_{l-1}) \right]-(t_{l}-t_{l-1})\left[W_0(t_{l+1})-W_0(t_{l})\right]\right\},$$
with 
$W_0$ standing for the standard Brownian bridge.
\item[(ii) -] Under the assumptions of Theorem \ref{th2},  under $H_{1,n}^k$, 
$\mathcal Z_n(t_1,t_2,\ldots,t_k)$ converges weakly to $\widetilde{\mathcal Z}_0(t_1,t_2,\ldots,t_k)$ with representation 
$$\widetilde{\mathcal Z}_0(t_1,t_2,\ldots,t_k) =
\mathcal Z_0(t_1,t_2,\ldots,t_k)+\sum_{l=1}^kA_l (t_{l+1}-t_l) (t_l -t_{l-1}).$$
\item[(iii) -] Under the assumptions of Theorem \ref{th3},  under $H_{1}^k$, 
$\mathcal Z_n(\cdot)$ diverges in probability to $\infty$.
\end{itemize}
%
% for
%any $s=(s_1, \ldots,s_k), \ t=(t_1, \ldots,t_k) \in \mathbb R^k$ 
%  by $\Gamma(s,t)=0$ if $s=(s_1, \ldots,s_k), \ t=(t_1, \ldots,t_k) \notin \Theta_k$ and for 
%$s=(s_1, \ldots,s_k), \ t=(t_1, \ldots,t_k) \in \Theta_k$,
%\begin{eqnarray}
%\Gamma(s,t)&=& \sigma \Big[\sum_{l=1}^k  (t_{l+1}-t_l)(s_{l+1}-s_l)(t_{l}\wedge s_{l}-t_{l}\wedge s_{l-1}-t_{l-1}\wedge s_{l} + t_{l-1}\wedge s_{l-1}) \nonumber\\
%&-& 
% \sum_{l=1}^k (t_{l+1}-t_l)(s_{l}-s_{l-1})
%(t_{l}\wedge s_{l+1}-t_{l}\wedge s_{l}-t_{l-1}\wedge s_{l+1}+t_{l-1}\wedge s_{l}) \nonumber\\
%&-& 
%\sum_{l=1}^k (t_{l}-t_{l-1}) (s_{l+1}-s_l)
%(t_{l+1}\wedge s_{l}-t_{l+1}\wedge s_{l-1}-t_{l}\wedge s_{l}+t_{l}\wedge s_{l-1}) \nonumber\\
%&+& 
%\sum_{l=1}^k (t_{l}-t_{l-1}) (s_l-s_{l-1})(t_{l+1}\wedge s_{l+1}-t_{l+1}\wedge s_{l}-t_{l}\wedge s_{l+1} + t_{l}\wedge s_{l}) \nonumber\\
%&+& 
%2\sum_{l=1}^k \sum_{j \le l} (t_{l+1}-t_l)(s_{j+1}-s_j)
%(t_{l}\wedge s_{j}-t_{l}\wedge s_{j-1}-t_{l-1}\wedge s_{j} + t_{l-1}\wedge s_{j-1})
%\nonumber\\
%&-& 
%2\sum_{l=1}^k \sum_{j \le l}  (t_{l+1}-t_l)(s_{j}-s_{j-1})
%(t_{l}\wedge s_{j+1}-t_{l}\wedge s_{j}-t_{l-1}\wedge s_{j+1}+t_{l-1}\wedge s_{j})
%\nonumber\\
%&-& 
%2\sum_{l=1}^k \sum_{j \le l} (t_{l}-t_{l-1})(s_{j+1}-s_{j})
%(t_{l+1}\wedge s_{j}-t_{l+1}\wedge s_{j-1}-t_{l}\wedge s_{j}+t_{l}\wedge s_{j-1})
%\nonumber\\
%&+& 
%2\sum_{l=1}^k \sum_{j \le l}  (t_{l}-t_{l-1})(s_{j}-s_{j-1})
%(t_{l+1}\wedge s_{j+1}-t_{l+1}\wedge s_{j}-t_{l}\wedge s_{j+1} + t_{l}\wedge s_{j}) \Big], 
%\end{eqnarray}
\end{thm}
\bigskip

\noindent
{\bf Proof} - 
 Theorem \ref{th1} holds for any kernel $h$. In the particular case where  $h$ is such that its associated 
$h_{\dot F_l,1}$ and $h_{\dot F_l,2}$ satisfy
$h_{\dot F_l,1}(x)=-h_{\dot F_l,2}(x)$, one has that $W_1(t)=-W_2(t)$.  Theorem \ref{th5} then follows by taking into account this equality in Theorem \ref{th4}.
\bigskip

\begin{rem}
The common covariance kernel of $\mathcal Z_0$ and $\widetilde Z_0$ is defined by $\Gamma(s,t), s, t \in \mathbb R^k$ with 
$ \sigma_{11}=\sigma_{22}=\sigma^2$ and 
$ \sigma_{12}=\sigma_{21}=-\sigma^2$. 
\end{rem}

\noindent
\begin{rem}
It is easy to check that anti-symmetric kernels $h$ are such that their associated $h_{\dot F_l,1}$ and $h_{\dot F_l,2}$ satisfy the property $h_{\dot F_l,1}(x)=-h_{\dot F_l,2}(x)$.
\end{rem}
%
%\begin{cor} \label{cor2}
%Assume that the assumptions of Theorem \ref{th2} hold, and that $h$ is such that its associated  $h_{F,1}$ and $h_{F,2}$ satisfy$h_{F,1}(x)=-h_{F,2}(x)$. Then 
 %
%\begin{itemize}
%\item[i- ] Under $H_0$, as $n$ tends to infinity, one has the following convergence in distribution 
%\begin{eqnarray}
%$$\mathcal T_{1,n} \longrightarrow  \sigma \sup_{\lambda \in [0,1]}\left | W^0(\lambda) \right |$$
%$$\mathcal T_{2,n} \longrightarrow \sigma^2 \sum_{j \ge 1} {1 \over j^2 \pi^2} \chi_j^2$$
%%\end{eqnarray}
%\item[ii- ] Under $H_{1,n}^k$, as $n$ tends to infinity, one has the following convergence in distribution 
%\begin{eqnarray}
%$$\mathcal T_{1,n} \longrightarrow \sup_{\lambda \in [0,1]}\left |(1-\lambda) \lambda A+ \sigma W^0(\lambda) \right |$$
%$$\mathcal T_{2,n} \longrightarrow \sum_{j \ge 1} {1 \over j^2 \pi^2} \chi_j^{*2},$$
%\end{eqnarray}
%\end{itemize}
%where $W^0$ is the Brownian bridge on $[0,1]$, the $ \chi_j^{2}$'s and $ \chi_j^{*2}$'s are as in Theorem \ref{th4} but the non-centrality parameters are  $2A^2 \left\{{2[1-(-1)^j] / j\pi} \right\}^2 \sigma^{-2}$.
%\end{cor}
%
%
%
%\begin{rem} \label{rm3} It is easy to check that anti-symmetric kernels $h$ are such %that their associated $h_{F,1}$ and $h_{F,2}(x)$ satisfy the property $h_{F,1}(x)=-h_{F,2}(x)$.
%\end{rem}
%
\medskip

 There are various techniques for estimating long-run variances. 
%Some use 
The subsampling approach introduced by \cite{carl} (see, e.g., \cite{schet}) is based on the  idea that, 
%Some others use the spectral density of the series under study. Concretely, the general
 %idea behind the subsampling methods is as follows, 
given a stationary time series $(\vartheta_i)$, its associated long-run variance $\sigma_\vartheta^2$ is defined as $\lim_{n\rightarrow \infty}\text{Var}(\sqrt{n} \overline \vartheta_n)$$=\sum_{j=-\infty}^{\infty} \gamma_\vartheta(j)=\gamma_\vartheta(0)+2\sum_{j=1}^{\infty} \gamma_\vartheta(j)$, where $\gamma_\vartheta$ is the autocorrelation function of $(\vartheta_i)$ and $\overline \vartheta_n$ is the sample mean of $\vartheta_1, \ldots, \vartheta_n$. Then, for some suitable integer $b_n$ ($b_n<n$), subsamples of lenghts, say $l_{n_j}<n$, $1\le j \le b_n$ are obtained from $\vartheta_1, \ldots, \vartheta_n$. From these sample, $b_n$ copies $\sqrt{l_{n_j}} \overline \vartheta_{l_{n_j}}$, $j=1, \ldots, b_n$ of $\sqrt{n} \overline \vartheta_n$ are computed. The sample variances of these copies is taken as an estimator of $\sigma_\vartheta^2$.
% is obtained as the sample mean of these sample variances. 
From the latter expression of $\sigma_\vartheta^2$, another class of estimators proposed by \cite{newey} is given as
$$\widehat \sigma_\vartheta=\widehat \gamma_\vartheta(0)+2\sum_{j=1}^{\infty} K\left({j \over \ell_n}\right)\widehat \gamma_\vartheta(j),$$
where $\ell_n$ is a bandwidth and $K$ is a symmetric kernel such that $K(0)=1$. The last approach that we present is as follows. Denoting the spectral density of $(\vartheta_i)$ by $f_\vartheta$, and recalling that $\sigma_\vartheta^2=2\pi f_\vartheta(0)$, an estimator of $\sigma_\vartheta^2$ is obtained by plugging an estimator of $f_\vartheta(0)$ in the preceding equality. A candidate estimator of $f_\vartheta(0)$ is 
$$\widehat f_\vartheta(0)={1 \over 2 \pi}\sum_{|j|\le n-1}w_n(j) \widehat \gamma_\vartheta(j),$$
where the $w_n(j)$ are weights such that $w_n(j) \longrightarrow 0$, as $|j| \rightarrow \infty$, and the  $\widehat \gamma_\vartheta(j)$ are the usual sample autocovariances. We can also note that if $\sigma_\vartheta^2$ is expressed as a function of unknown parameters whose estimators can be found, based on these estimators, an estimator of $\sigma_\vartheta^2$ can be obtained. This is the case for the $\vartheta_i$'s generated by an AR(1) model, as will be seen in the next section.
\medskip

The consistency of the above classes of estimators is established for a broad class of stationary time series. 
Let $\widehat \sigma^2$ be any such consistent estimator of $\sigma^2$ under $\mathcal H_0$, converging in probability to a positive number $\sigma^*$ under $\mathcal H_1^k$.
 Consider the statistics 
$$\widetilde{\mathcal T}_{1,n}={\mathcal T_{1,n} \over \widehat \sigma} \ \mbox{  and  } \ \widetilde{\mathcal T}_{2,n}={\mathcal T_{2,n} \over \widehat \sigma^2}.$$
%\medskip
%
%
\begin{cor} \label{cor4}
Assume that the assumptions of Theorem \ref{th2} hold. Let $(\mathcal B(s) : s\in [0,1]^k)$ be the process defined in Theorem \ref{th5}. Then 
\begin{itemize}
\item[i- ] Under $\mathcal H_0$, as $n$ tends to infinity, one has the following convergence in distribution, 
%\begin{eqnarray}
$$\widetilde{\mathcal T}_{1,n} \longrightarrow  \sup_{s \in [0,1]^k}\left | \mathcal B(s) \right | \ \mbox{ and } \ \widetilde{\mathcal T}_{2,n} \longrightarrow \int_{[0,1]^k} \mathcal B^2(s) ds = \sum_{j \ge 1} \widetilde\zeta_j \chi_j^2,$$
where the $ \chi_j^{2}$'s are iid chi-square random variables with one degree of freedom and the $\widetilde\zeta_j$'s are  the eigenvalues of the operator $\mathcal K_\Gamma/ \sigma^2$.
%
%\end{eqnarray}
\item[ii- ] Under $\mathcal H_{1,n}^k$, as $n$ tends to infinity, one has the following convergence in distribution, 
%\begin{eqnarray}
$$\widetilde{\mathcal T}_{1,n} \longrightarrow 
\sup_{s \in [0,1]^k}\left | \widetilde{\mathcal Z}_0(s) \right | \ \mbox{ and } \ \widetilde{\mathcal T}_{2,n} \longrightarrow  \int_{[0,1]^k} \widetilde{\mathcal Z}_0^2(s) ds = \sum_{j \ge 1} \widetilde\zeta_j \chi_j^{*2},$$
where the $ \chi_j^{*2}$'s are iid non-central chi-square random variables with one degree of freedom and non-centrality parameters $\rho_j^2\widetilde\zeta_j^{-1}\sigma^{-2}$
with the $\rho_j$'s defined in Theorem \ref{th4} and the $g_j$'s appearing in their expressions being the eigenfunctions of the integral operator $\mathcal K_\Gamma/\sigma^2$, associated with the eigen-value $\widetilde \zeta_j$.
%=\sum_{l=1}^kA_l\int_0^1 \ldots \int_0^1 (s_{l+1}-s_l) (s_l -s_{l-1})
%g_j(s_1, \ldots,s_k)\nbI((s_1,\ldots,s_k) \in \Theta_k)  ds_1 \ldots ds_k.$$
%\end{eqnarray}
\item[iii- ] Under $\mathcal H_1^k$, as $n$ tends to infinity, one has the following convergence in probability, 
%\begin{eqnarray}
$$\widetilde{\mathcal T}_{1,n} \longrightarrow \infty, \ \ \ \widetilde{\mathcal T}_{2,n} \longrightarrow \infty.$$
%\end{eqnarray}
\end{itemize}
\end{cor}
\bigskip

\noindent
{\bf Proof} -  Using the consistency of  $\widehat \sigma^2$ to $\sigma^2$ and  
 Theorem \ref{th5}, proceeding as in the proof of this theorem with $\mathcal Z_n/ \widehat \sigma$ substituted for $\mathcal Z_n$, one easily establishes this corollary.
\bigskip

\begin{rem} Using Theorem \ref{th4}, the asymptotic distributions of
%$\mathcal T_{1,n}$, $\widetilde{\mathcal T}_{1,n}$, 
$\mathcal T_{2,n}$ and $\widetilde{\mathcal T}_{2,n}$ under $\mathcal H_0$ and $\mathcal H_{1,n}^k$ can be approximated, even for more general kernel $h$.  This is an advantage over $\mathcal T_{1,n}$ and 
$\widetilde{\mathcal T}_{1,n}$ whose asymptotic distributions are more difficult to approximate. However, since $\widetilde{\mathcal T}_{1,n}$ and 
$\widetilde{\mathcal T}_{2,n}$ are asymptotic distribution-free under the null hypothesis, in the next section, we only consider kernels leading to these statistics. That is, the tests studied in that section are based on $\widetilde{\mathcal T}_{1,n}$ and 
$\widetilde{\mathcal T}_{2,n}$ whose asymptotic null distributions are free. Then, the quantiles of their limiting distributions are approximated by their simulated empirical counter part.
\end{rem}
\section{Simulation study}
For the simulations, we use the version 4.4.2 of the software R. 
As mentioned in Remark 6 above, we study only the tests based on $\widetilde{\mathcal T}_{1,n}$ and $\widetilde{\mathcal T}_{2,n}$. Potential kernels are anti-symmetric kernels, in particular $h(x,y)= \nbI(x<y)$ and $h(x,y)= x-y$. For the first, it is straightforward that 
 $\theta(F_1,F_1)=1/2$ and $h_{F_1,1}(x)=1/2-F_1(x)=-h_{F_1,2}(x)$. For the latter, it is a trivial matter that $\theta(F_1,F_1)=0$ and that  $h_{F_1,1}(x)=-h_{F_1,2}(x)$ as $h$ is anti-symmetric. The corresponding $\sigma^2$ are respectively 
$$ \sigma_1^2=\sigma^2= \mbox{Var}\left[  F_1(X_1) \right]+2 \sum_{j\ge1}\mathbb{C}\mbox{ov}(F_1(X_1), F_1(X_{1+j}))$$
and  
$$ \sigma_2^2=\sigma^2=  \mbox{Var}\left(X_1 \right)+2 \sum_{j\ge1}\mathbb{C}\mbox{ov}(X_1, X_{1+j}).$$
We restrict to $h(x,y)= x-y$ and $k=2$, the case $k=1$ being treated in \cite{nga}. Our results are applied to simple models for detecting : (i) two changes in the mean of a shifted white noise, (ii) two changes in the mean of an AR(1), (iii) one change in the mean and one change in the autocorrelation of an AR(1), and (iv) one change in the mean and one change in the variance of shifted white noise. 
%In other words, we check the difference between the distributions of the observations by checking the differences between their means. 
\bigskip

To obtain the critical values of the tests, we simulated 5000 observations of the stochastic process $\mathcal B(t_1,t_2)$ defined in Theorem \ref{th5} for $k=2$, and considered empirical quantiles of the limiting distribution as critical values. To simulate the Brownian bridge, following Donsker's theorem, we approximated the Brownian motion by the normalized partial sum process 
$$S_m(t)={1 \over \sqrt{m}} \sum_{i=1}^{[mt]} \eta_i, \ t \in [0,1], \ m=2000,$$
with the $\eta_i$'s simulated from a standard Gaussian distribution. With this, the Brownian bridge was approximated by 
$S_m(t)-tS_m(1), \ t \in [0,1]$. The process $\mathcal B (t_1,t_2)$ was then approximated by substituting $S_m(t)-tS_m(1)$ for $W_0(t)$. For computing the values of $\sup_{0\le t_1<t_2 \le 1} |\mathcal B(t_1,t_2)|$, we considered the couples $(i/m,j/m)$, $j=2, \ldots,m-2$, \ $i=1, \ldots, j-1$ and took it as the random variable \break $\max_{2\le j \le m-2, 1\le i \le j-1} |\mathcal B(i/m, j/m)|$. For computing the values of $\int_{\Theta_2} \mathcal B^2(t_1,t_2) dt_1 dt_2$, we approximated it by the random variable $\sum_{2\le j \le m-2, 1\le i \le j-1}\mathcal B^2(i/m, j/m)/m^2$. 
 Using all these approximations, and based on 5000 replications, the quantiles of $\sup_{0\le t_1<t_2 \le 1} \mathcal |B (t_1,t_2)|$ and $\int_{\Theta_2} B^2(t_1,t_2) dt_1 dt_2$ are approximated by the values in Table \ref{tab1}. 
%
%\bigskip
\begin{table}[ht]
\centering
\caption{Empirical quantiles of the limiting distributions of the test statistics.}
\label{tab1}
\begin{tabular}{|lr|cc|}
\noalign{\smallskip}\hline\noalign{\smallskip}
&Quantile&KS&CV\\
Level &&&\\ 
\hline 
%\smallskip
 {\bf 0.01}&&{\bf 1.66}&{\bf 0.249}\\
0.02&&1.55&0.206\\  
0.03 &&1.45&0.182\\ 
0.04&&1.43&0.162\\
{\bf 0.05}&&{\bf 1.38}&{\bf 0.145}\\
0.06&&1.34&0.134 \\  
0.07&&1.32&0.126\\
0.08&&1.29&0.119\\ 
0.09&&1.27&0.113\\ 
{\bf 0.10}&&{\bf 1.26}&{\bf 0.107}\\  
\hline
\end{tabular}
\end{table}

%\noindent
For evaluating the performances of our tests, we sampled 1000 sets of $n=200$ data $X_1, X_2, \ldots, X_n$ with $X_i=Y_{1,i}, i=25,\dots,100$,
$X_{75+i}=Y_{2,25+i}, i=1, \dots, 75$, $X_{150+i}=Y_{3,25+i}, i=1, \dots, 50$,  where 
 the $Y_{i,j}$'s are from the piece-wise stationary models
\begin{equation} \label{mod}
%X_i= \left \{
\begin{array}{cc}
Y_{1,j}=\mu_1+\rho_1 Y_{1,j-1}+\omega _1\varepsilon_j & j=1, \dots, 100,\\
Y_{2,j}=\mu_2+\rho_2 X_{j-1}+\omega _2\varepsilon_j & j=1, \dots, 100,\\
Y_{3,j}=\mu_3+\rho_3 X_{j-1}+\omega_3 \varepsilon_j & j=1, \dots, 100,
\end{array}
%\right.
\end{equation}
with the $\mu_j$'s, $\rho_j$'s and the $\omega_j$'s being real numbers, the $\varepsilon_i$'s iid random variables  and for all $i=1, \ldots, 200$, $\varepsilon_i \sim \mathcal N(0,1)$. These models are studied in \cite{yau} and their stationarity and ergodicity are studied in  \cite{nga1}.

We considered the following scenarios : 
\begin{itemize}
\item[(i)] Two changes in the mean : $\mu_1=0$, $\mu_2, \mu_3 \ne 0$,  $\mu_2 \ne \mu_3$, $\rho_i=0$ and $\omega_i=1$, $i=1, 2,3$.
\item[(ii)] Two changes in the mean : $\mu_1=0$, $\mu_2, \mu_3 \ne 0$, $\mu_2 \ne \mu_3$, $\rho_i=0.2$  and $\omega_i=1$, $i=1, 2,3$.
\item[(iii)] One change in the mean and one change in the autocorrelation : $\mu_1=0$, $\mu_2=\mu_3 \ne 0$, $\rho_1=\rho_2=0$, $\rho_3 \ne 0$ and $\omega_i=1$, $i=1, 2,3$.
\item[(iv)] One change in the mean and one change in the variance : $\mu_1=0$, $\mu_2=\mu_3 \ne 0$, $\rho_i=0$, $i=1, 2,3$,  
$\omega_1=\omega_2=1$ and $\omega_3 \ne 1$.
\end{itemize}

%\noindent
Here, the null model is an AR(1) of the form 
$$X_i=\mu+\rho X_{i-1}+\omega \varepsilon_i. $$
From simple computations, the corresponding long-run variance is 
$$\sigma^2={\omega^2 \over (1-\rho)^2}.$$
Let $\widehat \omega^2$ and $\widehat \rho$ be consistent estimators of $ \omega^2$ and $\rho$. 
% Examples of such estimators are the Yule-Walker estimators (see, e.g.,, \cite{bro}) given by 
An estimator $\widehat \sigma^2$ of $\sigma^2$ can be obtained by substituting these quantities in the above expression.
%$$\widehat \sigma^2={\widehat \omega^2 \over (1-\widehat \rho)^2},$$
%
%where $\widehat \omega^2$ and $\widehat \rho$ are consistent estimators of $ \omega^2$ and $\rho$.
 Examples of $ \widehat \omega^2$ and $\widehat \rho$ are the Yule-Walker estimators (see, eg, \cite{bro}) given by 
$$\left\{
\begin{array}{l} \displaystyle 
\widehat \rho = {\sum_{t=1}^{n-1}(X_{t-1}-\overline X)(X_t-\overline X) \over \sum_{j=1}^n(X_t-\overline X)^2}\\
\medskip

\widehat \mu = \overline X(1- \widehat \rho)\\
\smallskip

\displaystyle
\widehat \omega^2= {1 \over n} \sum_{t=1}^n (X_t-\widehat \mu- \widehat \rho X_{t-1})^2.
\end{array}
\right.$$

In what follows, we use the estimator of $\sigma^2$ obtained this way, because it is less subjective than the ones described in Subsection 3.2.
With it, the KS-type statistic $\mathcal T_{1,n}$ and the CV-type statistic $\mathcal T_{2,n}$ were computed by taking respectively, the maximum and the sum of the $X_i-X_j$'s over the couples $(i,j)$, $j=2, \ldots,n-2$, $i=1, \ldots, j-1$. Then the distribution-free test statistics $\widetilde{\mathcal T}_{1,n}$ and $\widetilde{\mathcal T}_{2,n}$ were obtained by dividing the resulting quantities, respectively by $\widehat \sigma$ and $\widehat \sigma^2$.

 For each of the 1000 samples generated from (\ref{mod}), the values of these distribution-free statistics were respectively compared to the 0.95-quantiles of the limiting distribution of the KS and CV tests. The ratios of the number of times they were larger than these numbers, over 1000 are their empirical powers. These were compared to those of the CUSUM test, the R\'enyi-type test and the Darling-Erd\"os test presented in \cite{laj3}, and  whose codes are available in the package CPAT of the software R. Table \ref{tab20} displays the empirical levels of the tests for some particular null models, while Tables \ref{tab2} and \ref{tab3} list their empirical powers for all the scenarios studied. From Table \ref{tab20}, it can be seen that the empirical level of the RC test ranges from 0.071 to 0.734. That of the DE test is nil for almost all the models except for the zero-mean AR(1) models  (lines 9-13). The empirical level of the remaining tests are more closer to the nominal level of 5\% except for the pure AR(1) (lines 7-13). While those of the CU test are too large besides the nominal level, those of the KS test are rather too small.

For changes in the mean (see Table \ref{tab2}), all the tests have nice powers, except the DE test for which the local power is generally too poor. For changes in the mean and correlation or variance (see Table \ref{tab3}), the power of the RC test is generally smaller than the others;  that of the DE test is unstable, while those of the remaining tests are nice and competitive. 
% One can also observe that for change detection in the mean, for $\rho_1=0$ (see the first two columns of Table \ref{tab2}) and $\rho_1=0.2$ (see the last two columns of Table \ref{tab2}), the power of the CV-test is larger against local alternatives and for jumps moving in the same direction. The power of the KS-test is smaller in the vinicity of the null hypothesis and is larger for jumps in the opposite directions. The same observations can be done for mixed changes (see Table \ref{tab3}).
%
%For changes in the mean, the power of both tests against detecting two changes in the autocorrelation or two changes in the variance was too poor. Thus, the tests seemed not to be able to detect these type of changes. The reason for this may be that the kernel $k(x,y)=x-y$ used in the test statistics is more adapted to changes in the location parameters, as, for example, the mean.
%
%
\begin{center}
%\footnotesize{
\begin{table}
\caption{Empirical level of the tests for scenarios from model (\ref{mod}) : $\mu_1=\mu_2=\mu_3$, $\rho_1=\rho_2=\rho_3$ and $\omega_1=\omega_2=\omega_3$ (no change neither in the mean nor in the autocorrelation and in the variance.}
\label{tab20}
\begin{tabular}{|l|ccccc|}
%\hline\noalign{\smallskip}
\noalign{\smallskip}\hline\noalign{\smallskip}
$(\mu_1, \rho_1,\omega_1)$ & KS&CV&CU&RC&DE
\\ 
\hline 
%\smallskip
$(0.1,0,1)$ &0.042 & 0.052&0.051&{\bf 0.098}&0.000\\  
$(0.4,0,1)$ & 0.047& 0.048&0.055&{\bf 0.074}&0.000\\ 
$(-0.8,0,1)$&  0.042& 0.046&0.041&{\bf 0.081}&0.000\\
$(1.2,0,1)$& 0.049&0.057&0.055&{\bf 0.109}&0.000\\
$(-1.6,0,1)$&0.046&0.052&0.053&{\bf 0.121}&0.000\\  
$(2,0,1)$ &0.046&0.052&0.043&{\bf 0.152}&0.000\\
$(0,0.1,1)$& 0.030& 0.035&0.035&{\bf 0.111}&0.003\\ 
$(0,0.3,1)$&0.032& 0.038&0.062&{\bf 0.156}&0.002\\ 
$(0,0.5,1)$ & {\bf 0.024} & 0.034&{\bf 0.139}&{\bf 0.244}&{\bf 0.018}\\  
$(0,0.7,1)$ & {\bf 0.010}&0.032&{\bf 0.287}&{\bf 0.483}&{\bf 0.087}\\ 
$(0,0.9,1)$&  {\bf 0.010}&0.039&{\bf 0.719}&{\bf 0.734}&{\bf 0.486}\\
$(0.7,0.4,1)$& {\bf 0.020}&0.032&{\bf 0.074}&{\bf 0.224}&0.002\\
$(-1,0.6,1)$&{\bf 0.019}& 0.029&{\bf 0.174}&{\bf 0.380}&0.031\\ 
%$(-0.90,0.90)$&0.870 & 0.917&0.901&0.143&0.179\\  
%$(0,0,0.1)$ &0.019& 0.028&0.175&0.380&0.032\\
$(0,0,0.2)$& 0.042& 0.052&0.050&{\bf 0.078}&0.000\\ 
$(0,0,0.4)$&0.040& 0.054&0.043&{\bf 0.083}&0.000\\ 
$(0,0,0.6)$&0.045& 0.047&0.042&{\bf 0.072}&0.000\\ 
$(0,0,08)$&0.050& 0.052&0.051&{\bf 0.082}&0.000\\ 
$(0,0,1)$&0.045& 0.049&0.047&{\bf 0.071}&0.000\\ 
$(0,0,1.5)$&0.043& 0.049&0.051&{\bf 0.081}&0.000\\ 
$(0,0,2)$&0.041& 0.045&0.049&{\bf 0.082}&0.001\\ 
\hline
\end{tabular}
\end{table}
%}
\end{center}

%
%
%{\footnotesize
\begin{center}
%\footnotesize{
\begin{table}
\caption{Power of the tests for scenarios from model (\ref{mod}) : $\mu_1=0$, $\mu_2, \mu_3 \ne 0$,  $\mu_2 \ne \mu_3$, $\rho_i=0$ and $\omega_i=1$, $i=1, 2,3$ (left side of the table) ; $\mu_1=0$, $\mu_2, \mu_3 \ne 0$, $\mu_2 \ne \mu_3$, $\rho_i=0.2$  and $\omega_i=1$, $i=1, 2,3$ (right side of the table).}
\label{tab2}
\begin{tabular}{|l|ccccc|ccccc|}
%\hline\noalign{\smallskip}
\noalign{\smallskip}\hline\noalign{\smallskip}
$(\mu_2,\mu_3)$ & KS&CV&CU&RC&DE&KS&CV&CU&RC&DE
\\ 
\hline 
%\smallskip
$(0.00,0.00)$ &0.042 & 0.051&0.048&0.091&0.002&0.039&0.052&0.051&0.133& 0.000\\  
$(0.01,-0.05)$ & 0.046& 0.053&0.043&0.071&0.000&0.038&0.045&0.050&0.129&0.000\\ 
$(-0.03,0.05)$&  0.049& 0.056&0.051&0.081&0.000&0.037&0.046&0.047&0.107&0.000\\
$(0.06,-0.05)$& 0.050&0.058&0.054&0.077&0.000&0.034&0.047&0.067&0.135&0.001\\
$(0.06,-0.07)$&0.056&0.058&0.060&0.082&0.000&0.037&0.048&0.065&0.135& 0.001\\  
$(0.08,-0.09)$ &0.066&0.063&0.072&0.085&0.001&0.042&0.058&0.064&0.128&0.001\\
$(0.10,0.14)$& 0.080& 0.105&0.089&0.086&0.001&0.074&0.100&0.128&0.135&0.001\\ 
$(0.15,0.20)$&0.158& 0.193&0.178&0.083&0.001&0.096&0.132&0.135&0.144&0.002\\ 
$(-0.50,0.5)$ &0.854 & 0.376&0.945&0.203&0.444&0.722&0.274&0.948&0.229&0.358\\  
$(0.25,0.20)$ &  0.057&0.057&0.054&0.081&0.000&0.158&0.183&0.253&0.152&0.003\\ 
$(0.50,-0.30)$&  0.223&0.230&0.245&0.093&0.004&0.365&0.088&0.689&0.175&0.079\\
$(-0.80,0.40)$& 0.549&0.132&0.712&0.106&0.081&0.747&0.095&0.986&0.223&0.489\\
$(-0.45,-0.60)$&0.870& 0.182&0.987&0.223&0.482&0.786&0.876&0.920&0.195&0.197\\ 
$(-0.90,0.90)$&0.870 & 0.917&0.901&0.143&0.179&1.000&0.770&1.000&0.916&0.998\\  
$(1.00,-0.80)$ &1.000& 0.840&1.000&0.928&0.998&1.000&0.560&1.000&0.840&0.995\\
$(1.00,1.50)$& 1.000& 0.659&1.000&0.874&0.994& 1.000& 1.000&1.000&0.875&1.000\\ 
$(-1.50,1.00)$&1.000& 1.000&1.000&0.882&0.999& 1.000& 1.000&0.356&1.000&1.000\\ 
\hline
\end{tabular}
\end{table}
%}
\end{center}

\begin{table}
\caption{Power of the tests for scenarios from model (\ref{mod}) : $\mu_1=0$, $\mu_2=\mu_3 \ne 0$, $\rho_1=\rho_2=0$, $\rho_3 \ne 0$ and $\omega_i=1$, $i=1, 2,3$ (second and third columns); $\mu_1=0$, $\mu_2=\mu_3 \ne 0$, $\rho_i=0$, $i=1, 2,3$,  
$\omega_1=\omega_2=1$ and $\omega_3 \ne 1$ (last two columns).}
\label{tab3}
\begin{tabular}{|l|ccccc||l|ccccc|}
%\hline\noalign{\smallskip}
\noalign{\smallskip}\hline\noalign{\smallskip}
$(\mu_2,\rho_3)$ &KS&CV&CU&RC&DE&$(\mu_2,\omega_3)$ &KS&CV&CU&RC&DE
\\ 
\hline 
%\smallskip
(0.0,0.0)&0.047&0.052&0.047&0.075&0.000&(0.0,0.0)&0.046&0.056&0.044&0.083&0.001\\  
(0.1,0.1)&0.062&0.075&0.073&0.083&0.001&(1.0,0.5)&0.097&0.089&0.108&0.106&0.000\\ 
(-0.1,-0.3)&0.109&0.126&0.123&0.094&0.002&(-0.2,0.5)&0.238&0.230&0.269&0.114&0.011\\
(0.1,-0.5)&0.096&0.114&0.112&0.088&0.001&(-0.3,0.5)&0.469&0.476&0.509&0.124&0.022\\
(-0.1,0.7)&0.038&0.042&0.036&0.085&0.001&(0.5,0.5)&0.899&0.914 &0.923&0.169&0.277\\  
(0.1,0.9)&0.036&0.043&0.032&0.102&0.001&(0.1,1.5)&0.076&0.110&0.082&0.106&0.006 \\
(-0.5,-0.2)&0.779&0.797&0.836&0.111&0.129&(-0.3,1.5)&0.262&0.327&0.293&0.111&0.008\\ 
(0.5,0.4)&0.850&0.868&0.893&0.139&0.108&(0.5,1.5)&0.643&0.697&0.697&0.013&0.066\\ 
(-0.5,0.6)&0.787&0.817&0.834&0.122&0.078&(0.8,1.5)&0.990&0.972&0.991&0.208& 0.587\\  
(0.5,-0.8)&0.716&0.746&0.756&0.145&0.127&(1.0,1.5)&0.999&0.998&1.000&0.319&0.908\\ 
(-0.5,-0.9)&0.682&0.738&0.740&0.119&0.101&(-0.5,-1.5)&0.649&0.678&0.701&0.127&0.053\\
(1.0,0.1)&1.000&1.000&1.000&0.495&0.987&(0.4,-0.5)&0.701&0.709&0.731&0.142&0.100\\
(-1.0,-0.5)&1.000&1.000&1.000&0.482&0.957&(-0.7,-0.5)&0.997&0.998&1.000&0.253&0.777\\ 
(0.9,-0.7)&0.954&0.968&0.964&0.159&0.219&(1.0,-2.0)&0.998&0.994&0.998&0.223&0.707\\  
(1.0,-0.9)&0.999&0.999&0.999&0.307&0.857&(-1.0,-0.2)&1.000&1.000&1.000&0.782&0.999\\
(-0.6,0.7)&1.000&0.997&1.000&0.376&0.924&(-0.2,0.1)&0.246&0.240&0.285&0.119&0.004\\ 
(0.8,0.8)&1.000&0.999&0.999&0.201&0.714&(0.9,2.0)&0.979&0.967&0.985&0.200&0.492\\ 
\hline
\end{tabular}
\end{table}

\section{Concluding remarks}
The current work is an extension of \cite{nga} to multiple change-points detection. 
The results are applied to testing for two changes of possibly different natures in a series simulated from a piece-wise stationary AR(1) model. The critical values were obtained as simulated empirical quantiles of the limiting distributions. For the kernel $h(x,y)=x-y$ that we used, we observed that our tests have nice empirical levels, and that their empirical powers compared to the CUSUM tests, the R\'enyi-type test and the Darling-Erd\"os test were very competitive. 

The powers of all the tests were generally too poor for the series with the two changes either in the autocorrelation or in the variance.  For our tests, this can be explained by the fact that the kernel used is more adapted to testing changes in the location parameters. This argument can also be invoked for the other tests.
%, it may come from the fact that they are built for testing 
% A more suitable kernel for change detection in the variance is $h(x,y)=(x-y)^2/2$. We did not try it because its use involves a lot of computations. 
%Also, we have not tried to implement the tests associated with more general kernels, because they seem to have a large cost computational time. All these remaining works will be the subject of further investigations.
%
\section{Appendix : Proofs of the results}
\subsection{Some useful results}
We need the following result proved by \cite{oo}.
\smallskip

Let $\xi_1,\xi_2,\ldots,\xi_n,\ldots$ be a strictly stationary sequence of zero-mean random variables, and let
$$
\sigma _{\ast}^{2}=\mathbb E(\xi_1^{2})+2\sum_{i=1}^{\infty} \mathbb E \left(\xi_1\xi_{i+1} \right).
$$

\begin{prop} \label{pro2}
Assume $ \mathbb E \left(\left\vert \xi_{i}\right\vert ^{2+\delta } \right)<\infty $ for some positive $\delta$ and 
$\xi_1,\xi_2,\ldots,\xi_n,\ldots$  is $\alpha$-mixing  with $\alpha$-rate satisfying 
\begin{equation*}
\sum_{i=1}^{\infty }[\alpha (i)]^{\frac{\delta }{2+\delta }}<\infty.
\end{equation*}
Then $\sigma _{\ast }^{2}<\infty.$ 
\medskip

\noindent
If  $\sigma_{\ast}>0$, then the sequence of processes
 $$S_{n}(t)=\frac{1}{\sigma_{\ast}\sqrt{n}}\sum_{i=1}^{[nt]}\xi_{i}, \ t \in [0,1]$$
converges weakly to a Wiener measure on $(D,\mathcal{D)}$, where 
$\mathcal{D}$ is the $\sigma $-fields of Borel sets for the Skorohod topology.
\end{prop}
\bigskip

\begin{lem} \label{lem3} (\cite{har1}) Let $\{\mathcal X_{ni}^*\}$ be a sequence of zero-mean absolutely regular
random variables (rv)'s with rates satisfying
\begin{equation} \label{hp1}
\sum_{n \geq 1}\big[\beta(n)\big]^{\delta/(2+\delta)} <\infty \  \mbox{ for some } \ \delta> 0.
\end{equation}
Suppose that for any $\kappa$, there exists a sequence $\{\mathcal Y^\kappa_{ni}\}$ of
rv's satisfying (\ref{hp1}) such that
\begin{equation} \label{hp2}
\sup_{n\in\mathbb{N}}\max_{0\leq i\leq n}|\mathcal Y^\kappa_{ni}|\leq B_\kappa < \infty,
\end{equation}
where $B_\kappa$ is some positive constant
\begin{equation} \label{hp3}
\sup_{n\in\mathbb{N}}\max_{0\leq i\leq n} \mathbb E\big(|\mathcal X^*_{ni}-\mathcal Y^\kappa_{ni}|^{2+\delta}\big)\longrightarrow 0  \
\mbox{ as } \ \kappa \rightarrow\infty
\end{equation}
\begin{equation} \label{hp4}
\frac{1}{n} \mathbb E\bigg[\Big(\sum^n_{i=1}\mathcal X^*_{ni}\Big)^2\bigg]\longrightarrow c \  \mbox{ as } \ n\rightarrow \infty,
\end{equation}
where c is some positive constant
\begin{equation} \label{hp5}
\frac{1}{n} \mathbb E\bigg[\Big(\sum^n_{i=1}\mathcal Y^\kappa_{ni}-\mathbb
E(\mathcal Y^\kappa_{ni})\Big)^2\bigg]\longrightarrow c_\kappa \  \mbox{ as } \ n\rightarrow \infty,
\end{equation}
where $c_\kappa$ is some constant $> 0$
\begin{equation} \label{hp6}
c_\kappa \longrightarrow c \ \mbox{ as } \ \kappa \rightarrow\infty.
\end{equation}
Then
\begin{equation*}
\frac{1}{\sqrt{n}}\sum^n_{i=1}\mathcal X^*_{ni}
\end{equation*}
converges in distribution to the normal distribution with mean 0 and variance $c$.
\end{lem}
 \bigskip

\begin{lem} \label{phil} (\cite{phil}) Probability measures on a product space are
tight iff all the marginal probability measures are tight on the component spaces.
\end{lem}
\subsection{Preliminaries}
In this subsection, we prove some preliminary results necessary to the proofs of
Theorems \ref{th1} and \ref{th2}.\bigskip

\begin{prop} \label{pro3}
Under the conditions of Theorem \ref{th1}, as $n$ tends to infinity, in probability, we have 
\begin{equation*}
n^{-3/2}\sup_{(t_1,t_2,\ldots,t_k) \in \Theta_k}\left|
\sum_{l=1}^k\sum_{i=[nt_{l-1}]+1}^{[nt_l]}\sum_{j=[nt_l]+1}^{[nt_{l+1}]}g_{F_1,F_1}(X_{i},X_{j})\right|
\longrightarrow 0.
\end{equation*}
Under the conditions of Theorem \ref{th2},  as $n$ tends to infinity, in probability, we have 
\begin{equation*}
n^{-3/2}\sup_{(t_1,t_2,\ldots,t_k) \in \Theta_k}\left|
\sum_{l=1}^k\sum_{i=[nt_{l-1}]+1}^{[nt_l]}
\sum_{j=[nt_l]+1}^{[nt_{l+1}]}g_{F^{(n)}_{l-1},F^{(n)}_l}(Y^{l-1}_{ni},Y^l_{nj})\right|
\longrightarrow 0
\end{equation*}
where by convention $[nt_0]=1$ and $[nt_{k+1}]=n$.
\end{prop}
\bigskip

\noindent
{\bf Proof} -
We only prove the first part. This needs two lemmas that we first state and prove.\bigskip

\begin{lem} \label{lem1}
Under the conditions of Theorem \ref{th1}, there exists a constant $Cst>0$ such that for any $(t_1,t_2,\ldots,t_k) \in \Theta_k$, 
\begin{eqnarray*}
&&\mathbb{E}\left\{\left[\sum_{l=1}^k\sum_{i=[nt_{l-1}]+1}^{[nt_l]}\sum_{j=[nt_l]+1}^{[nt_{l+1}]}
g_{F_1,F_1}(X_{i},X_{j})\right]^2 \right\}\\
&& \leq  Cst\Big\{\sum_{l=1}^k([nt_l]-[nt_{l+1}])([nt_{l+1}]-[nt_l])\\
&  & +2\sum_{1\leq l<l'\leq k} ([nt_l]-[nt_{l+1}])([nt_{l'}]-[nt_{l'-1}])\Big\},
\end{eqnarray*}
with the convention that  convention $[nt_0]=1$ and $[nt_{k+1}]=n$.
\end{lem}
\bigskip

\noindent
%{\bf Proof of Lemma \ref{lem1}} 
{\bf Proof} - For any $(t_1,t_2,\ldots,t_k) \in \Theta_k$, we can write
\begin{eqnarray*}
&&\mathbb{E}\left\{\left[\sum_{l=1}^k\sum_{i=[nt_{l-1}]+1}^{[nt_l]}\sum_{j=[nt_l]+1}^{[nt_{l+1}]}
g_{F_1,F_1}(X_{i},X_{j})\right]^2 \right\}\\ && \leq   \sum_{l=1}^k\sum_{i=[nt_{l-1}]+1}^{[nt_l]}\sum_{j=[nt_l]+1}^{[nt_{l+1}]}
\mathbb{E}\big\{\big[g_{F_1,F_1}(X_{i},X_{j})\big]^2\big\}\\
&  & + 2\sum_{1\leq l<l'\leq k}\sum_{1\leq i_1<i_2\leq [nt_l]}\sum_{[nt_{l'}]+1\leq j_1<j_2\leq[nt_{l'+1}]} H((i_1,j_1),(i_2,j_2))\\
& = & \sum_{l=1}^kA_n(t_l)+2\sum_{1\leq l<l'\leq k}B_n(t_l,t_{l'}),
\end{eqnarray*}
where
\begin{eqnarray*}
H((i_1,j_1),(i_2,j_2)) & = & \mathbb{E}\big\{[g_{F_1,F_1}(X_{i_1},X_{j_1})-h_{F_1,1}(X_{i_1})-h_{F_1,2}(X_{j_1})+\theta_h(F_1,F_1)]\\
&  & \times [g_{F_1,F_1}(X_{i_2},X_{j_2})-h_{F_1,1}(X_{i_2})-h_{F_1,2}(X_{j_2})+\theta_h(F_1,F_1)]\big\}.
\end{eqnarray*}
From the integrability condition, we have
$$\sup_{i,j\in\mathbb{N}}\mathbb{E}\big\{\big[g_{F_1,F_1}(X_1,X_2)\big]^2\big\}\leq Cst.$$
Then, for any $1 \le l \le k$, 
$$A_n(t_l)\leq Cst([nt_l]-[nt_{l-1}])([nt_{l+1}]-[nt_l]).$$
Similarly as in \cite{nga}, we prove that
$$B_n(t_l,t_{l'})\leq Cst ([nt_l]-[nt_{l+1}])([nt_{l'}]-[nt_{l'-1}])$$
and Lemma \ref{lem1} is proved.
\medskip

\noindent
For any $(t_1,t_2,\ldots,t_k) \in \Theta_k$, we now define the process $\mathcal{G}_n(t_1,\ldots,t_k)$ by
\begin{equation*}
\mathcal{G}_n(t_1,\ldots,t_k)=n^{-3/2}\sum_{l=1}^k\sum_{i=[nt_{l-1}]+1}^{[nt_l]}\sum_{j=[nt_l]+1}^{[nt_{l+1}]}
g_{F_1,F_1}(X_{i},X_{j}).
\end{equation*}
\begin{lem} \label{lem2}
Under the conditions of Theorem \ref{th1}, we have
\begin{equation} \label{E1}
\mathbb{E}\left[\left(\mathcal{G}_n(t_1^1,\ldots,t_1^k)-\mathcal{G}_n(t_2^{1},\ldots,t_2^k)
\right)^2 \right]\leq \frac{Cst}{n}\sum_{l=1}^k(t_1^l-t_2^l),
\end{equation}
for all $0< t_1^1<t_2^1<t_1^2<t_2^2<\ldots<t_1^k<t_2^k < 1$.
\end{lem}
\bigskip

\noindent
%{\bf Proof of Lemma \ref{lem2}}
{\bf Proof} - For any $0< t_1^1<t_2^1<t_1^2<t_2^2<\ldots<t_1^k<t_2^k < 1$, we can write
\begin{eqnarray*}
&&\mathbb{E}\left[\left(\mathcal{G}_n(t_1^1,\ldots,t_1^k)-\mathcal{G}_n(t_2^{1},\ldots,t_2^k)\right)^2 \right]\\
&& =  \mathbb{E}\bigg[\bigg(n^{-3/2}\sum_{l=1}^k\sum_{i=[nt_1^{l-1}]+1}^{[t_1^l]}\sum_{j=[nt_l]+1}^{[nt_1^{l+1}]}
g_{F_1,F_1}(X_{i},X_{j})\\
&  & - n^{-3/2}\sum_{l=1}^k\sum_{i=[nt_2^{l-1}]+1}^{[nt_2^l]}\sum_{j=[nt_2^l]+1}^{[nt_2^{l+1}]}
g_{F_1,F_1}(X_{i},X_{j})\bigg)^2 \bigg]\\
& \leq & 2n^{-3}\mathbb{E}\bigg[\bigg(\sum_{l=1}^k\sum_{i=[nt_2^{l-1}]+1}^{[nt_1^l]}\sum_{j=[nt_1^l]+1}^{[nt_2^{l}]}
g_{F_1,F_1}(X_{i},X_{j})\bigg)^2 \bigg]\\
&  & +2n^{-3}\mathbb{E}\bigg[\bigg(\sum_{l=1}^k\sum_{i=[nt_1^{l-1}]+1}^{[nt_2^l]}\sum_{j=[nt_2^{l}]+1}^{[nt_1^{l+1}]}
g_{F_1,F_1}(X_{i},X_{j})\bigg)^2\bigg].
\end{eqnarray*}
From Lemma \ref{lem1}, we deduce (\ref{E1}) and Lemma \ref{lem2} is proved.\bigskip

\noindent
For the simplification of the proof, we fix $k=2$. The case $k>2$ is lengthy be easy.
\medskip

\noindent
From Lemma \ref{lem2}, we deduce that for all $\epsilon>0$ and all $0\leq t^1_1 < t^1_2 < t_1^2 < t_2^2\leq1$,
$$\mathbb{P}\left(|\mathcal{G}_n(t^1_1, t^1_2)- \mathcal{G}_n(t^2_1, t^2_2)|\geq\epsilon \right)\leq {Cst\over \epsilon^2n}[(t_2^1- t^1_1)+(t_2^2- t_1^2)].$$
It implies that for all $0\leq s_1^1 < s_2^1 < s_1^2 < s_2^2\leq n$ with $0\leq s_1^1 \ll s_2^1 \ll s_1^2 \ll s_2^2\leq n$
\begin{eqnarray*}
\mathbb{P}\left[ \left|\mathcal{G}_n \left({s_1^1 \over n}, {s^1_2\over n}\right)- \mathcal{G}_n\left({s^2_1\over n}, {s^2_2\over n}\right) \right|\geq \epsilon\right] & \leq & {Cst \over \epsilon^2n^2}[(s_2^1- s^1_1)+( s_2^2- s_1^2)]\\
& \leq & {Cst \over \epsilon^2} \left[{(s_2^1- s_1^1)+(s_2^2- s_1^2) \over n^{5/4}} \right ]^{4/3}.
\end{eqnarray*}
Now consider the partial sum process defined by

$$\mathcal S_{s,t}=\sum_{2\leq j\leq t}\sum_{1\leq i\leq \min\{j-1,s\}}\mathcal A(i,j),$$
where
$$\mathcal A(i,j)= \displaystyle \left \lbrace
\begin{array}{lr}
\mathcal{G}_n\left({i\over n}, {j\over n}\right)- \mathcal{G}_n\left({i-1\over n}, {j\over n} \right)- \mathcal{G}_n\left({i\over n}, {j-1\over n}\right)+ \mathcal{G}_n\left({i-1\over n}, {j-1\over n}\right)&\mbox{ if } \  i<j \\
\mathcal{G}_n\left({i\over n}, {j\over n}\right)= 0 &\mbox{ if } \ i\geq j.
\end{array}
\right.$$
It results that
$$\mathcal S_{s,t}= \mathcal{G}_n\left({s\over n}, {t\over n}\right).$$
The last inequality is equivalent to
\begin{eqnarray*}
\mathbb{P}(|\mathcal S_{s,t} - \mathcal S_{s',t'}|\geq \epsilon) & \leq & {Cst \over \epsilon^2}\left\{{[(s'-s)+( t'-t)]\over n^{5/4}} \right\}^{4/3}\\
& \leq & {Cst \over \epsilon^2n^{1/3}}.
\end{eqnarray*}
From Theorem 10.2 of \cite{billingsley}, we easily deduce
$$\mathbb{P}\left(\max_{1\leq s <t\leq n-1}|\mathcal S_{s,t}|\geq \epsilon\right)\leq {Cst \over \epsilon^2n^{1/3}},$$
which implies that, as $n$ tends to infinity, in probability
\begin{equation*}
n^{-3/2}\sup_{0< t_1<t_2 < 1}\left|
\sum_{l=1}^2\sum_{i=[nt_{l-1}]+1}^{[nt_l]}\sum_{j=[nt_l]+1}^{[nt_{l+1}]}
g_{F_1,F_1}(X_{i},X_{j})\right|
\longrightarrow 0.
\end{equation*}
This completes the proof of Proposition \ref{pro3}.
\bigskip
%\end{prop}

\begin{prop} \label{pro2}
Under the conditions of Theorem \ref{th1}, as $n$ tends to infinity, we have the weak convergence result
\begin{equation} \label{p1}
\frac{1}{\sqrt n}  \sum_{i=1}^{[nt]}\left(
  \begin{array}{c}
   h_{F_1,1}(X_{i}) \\
   h_{F_1,2}(X_{i}) 
  \end{array}
\right)
\longrightarrow 
\left(
  \begin{array}{c}
   W_{1}(t) \\
   W_{2}(t) 
  \end{array}
\right),
\end{equation}
in the space $(D[0,1])^2$.
%where by convention $[nt_0]=1$ and $[nt_{k+1}]=n$.
\medskip

\noindent
Under the conditions of Theorem \ref{th2}, for any $l=1, \ldots, k$, as $n$ tends to infinity, we have the weak convergence result
\begin{equation} \label{p2}
\frac{1}{\sqrt n} \sum_{i=1}^{[nt]}\left(
  \begin{array}{c}
    h_{F^{(n)}_l,1}(Y^l_{ni}) \\
    h_{F^{(n)}_l,2}(Y^{l+1}_{ni}) \\
  \end{array}
\right)
\longrightarrow
\left(
  \begin{array}{c}
   W_1(t) \\
   W_2(t) 
  \end{array}
\right),
\end{equation}
in the space $(D[0,1])^2$.
%where by convention $[nt_0]=1$ and $[nt_{k+1}]=n$.
\end{prop}
\bigskip

\noindent
{\bf Proof} - 
We start with the proof of (\ref{p1}). Under the assumptions of Theorem \ref{th2}, the sequences $(h_{F_1,j}(X_i))_{i \in \mathbb Z}$, $j=1,2$ and  $\left( \left(
  \begin{array}{c}
   h_{F_1,1}(X_i) \\
   h_{F_1,2}(X_i) 
  \end{array}
\right)\right)_{i \in \mathbb Z}$ are strictly stationary. For any $a_1, a_2 \in \mathbb R$, one has by Proposition 1,
$$
\left \{
\begin{array}{l} \displaystyle 
\frac{1}{\sqrt n} \sum_{i=1}^{[nt]} h_{F_1,1}(X_i) \longrightarrow \sigma_{11} W_1(t)\\
\displaystyle
\frac{1}{\sqrt n} \sum_{i=1}^{[nt]} h_{F_1,2}(X_i) \longrightarrow \sigma_{22} W_2(t)\\
\displaystyle
\frac{1}{\sqrt n} \sum_{i=1}^{[nt]} \left[a_1 h_{F_1,1}(X_i) + a_2h_{F_1,2}(X_i)  \right] \longrightarrow \widetilde \sigma \widetilde W(t),
\end{array}
\right.$$
where $W_1, W_2$ and $\widetilde W$ are Wiener processes and $\widetilde \sigma$ equals $\sigma_{rp}$ with 
$a_1 h_{F_1,1}(X_i) + a_2h_{F_1,2}(X_i)$ substituted for $h_{F_1,1}(X_i)$ and $h_{F_1,2}(X_i)$. 
\smallskip

It is easy to see that in distribution $\widetilde \sigma \widetilde W(t)$ equals $ \sigma_{11} W_1(t)+ \sigma_{22} W_2(t)$, which establishes  (\ref{p1}).\medskip

For establishing (\ref{p2}) we cannot proceed directly as for (\ref{p1}) since the bidimensional series 
$\left( \left(\begin{array}{c}
    h_{F^{(n)}_l,1}(Y^l_{ni}) \\
    h_{F^{(n)}_l,2}(Y^{l+1}_{ni}) \\
  \end{array}
\right)\right)_{i \in \mathbb Z}$ is not strictly stationary. We need to study the tightness and finite-dimensional distributions of the bi-dimensional partial sum in the left-hand side of (\ref{p2}). 

For the study of the finite-dimensional distributions, by the Cram\'{e}r-Wold device it
suffices to show that for any $\ell \in \mathbb{N}^{\ast}$, any
$a_j, b_j, t_j\in \mathbb{R}$, $0=t_0<t_1<\ldots<t_\ell=1$
\begin{equation*}
\sum_{j=1}^{\ell}\frac{1}{\sqrt{n}}\sum_{i=[nt_{j-1}]+1}^{[nt_{j}]}\left[a_{j}h_{F_l^{(n)},1}(Y_{ni}^l)+b_{j}h_{F_l^{(n)},2}(Y_{ni}^{l+1}) \right]
\end{equation*}
converges in distribution to a Gaussian random variable.
\bigskip

For the simplification of the presentation of the proof,  we only treat the case $\ell=2$ and $0=t_0<t_1<t_2=1$. 
%$a_1<a_2$, $b_1<b_2$.
\\\\
The assumption (\ref{hp1}) readily holds from (\ref{bet}).
\medskip

Define, for $j=1,2$,
\begin{equation*}
\vartheta_{ni}^{(j)}=a_jh_{F_l^{(n)},1}(Y_{ni}^l)+b_jh_{F_l^{(n)},2}(Y_{ni}^{l+1}).
\end{equation*}
For establishing (\ref{hp4}), we need proving that, as $n$ tends to infinity,
\begin{equation*}
\mathbb E \left \{ \left[\frac{1}{\sqrt{n}}\left (\sum_{i=1}^{[nt_1]}
\vartheta_{ni}^{(1)}
+\sum_{i=[nt_1]+1}^{n} \vartheta_{ni}^{(2)}\right)\right]^{2} \right\}
% \longrightarrow_{n\rightarrow\infty }c
\end{equation*}
tends to some positive constant $c$.
\bigskip

\noindent
We have
\begin{eqnarray*}
&&\mathbb E\left\{\left[ \frac{1}{\sqrt{n}}\left(\sum_{i=1}^{[nt_1]}\vartheta_{ni}^{(1)}+\sum_{i=[nt_1]+1}^{n}\vartheta_{ni}^{(2)}\right) \right]^{2} \right\}
= \\
&  & \ \ 
\frac{1}{n}\left\{\mathbb E\left[\left(\sum_{i=1}^{[nt_1]}\vartheta_{ni}^{(1)}\right)^{2} \right]
+2\mathbb E\left[\left(\sum_{i=1}^{[nt_1]}
\vartheta_{ni}^{(1)}\right)\left(\sum_{i=[nt_1]+1}^{n}\vartheta_{ni}^{(2)}\right)\right]
+\mathbb E\left[ \left(\sum_{i=[nt_1]+1}^{n}\vartheta_{ni}^{(2)}\right)^{2}\right]\right\}.
\end{eqnarray*}
Since the random variables $\vartheta_{ni}^{(1)}$ and $\vartheta_{ni}^{(2)}$ are centered, we obtain
\begin{equation*}
\mathbb E\left[\bigg(\sum_{i=1}^{[nt_1]}\vartheta_{ni}^{(1)}\bigg)^{2} \right]
=[nt_1] \mathbb E\left[\big(\vartheta_{n1}^{(1)}\big)^{2} \right]+2\sum_{i=1}^{[nt_1]}
\sum_{j=1}^{[nt_1]-i}\mathbb E\left(\vartheta_{ni}^{(1)}\vartheta_{n,i+j}^{(1)}\right).
\end{equation*}

\noindent
From the condition of Theorem \ref{th2}, we deduce that
$\mathbb E\left[ \left(\vartheta_{n1}^{(1)}\right)^{2+\delta}\right]<\infty$, which implies that
\begin{equation*}
\sup_{n, i,j\geq 1}\left\vert \mathbb E\left(\vartheta _{ni}^{(1)}\vartheta_{n,i+j}^{(1)}\right)\right\vert  \leq
\beta^{\frac{\delta}{2+\delta}}(j)\left\{\mathbb E\left(\left|\vartheta_{ni}^{(1)}\right|^{2+\delta}\right)\right\}^{\frac{1}{2+\delta}}
\left\{\mathbb E\left(\left|\vartheta_{n,i+j}^{(1)}\right|^{2+\delta}\right)\right\}^{\frac{1}{2+\delta}}.
\end{equation*}
We get
\begin{equation*}
\mathbb E\left[\left(\sum_{i=1}^{[nt_1]}\vartheta _{ni}^{(1)}\right)^{2} \right]\leq \lbrack nt_1]
\mathbb E\left [ \left(\vartheta_{n1}^{(1)}\right)^{2} \right]+2[nt_1]\sum_{j=1}^{[nt_1]}\beta^{\frac{\delta}{2+\delta}}(j)M^{2},
\end{equation*}
where $M=\sup_{n \ge 1}\left\{\mathbb E\left[\left(\vartheta_{n1}^{(1)}\right)^{2+\delta}\right]\right\}^{\frac{1}{2+\delta}}$.
\\\\
It results that
\begin{eqnarray} \label{chp1}
&&\lim_{n\rightarrow \infty}\frac{1}{n}\bigg\{[nt_1]\mathbb E\left[\left(\vartheta_{n1}^{(1)}\right)^{2} \right]
+2\sum_{i=1}^{[nt_1]}\sum_{j=1}^{[n t_1]-i}\Big|\mathbb E\left(\vartheta_{ni}^{(1)}\vartheta_{n,i+j}^{(1)}\right)\Big|\bigg\} \nonumber \\
&\leq& t_1\Big\{\mathbb E\left[\left(\vartheta_{n1}^{(1)}\right)^{2}\right]+2\sum_{j=1}^{\infty}\beta^{\frac{\delta}{2+\delta}}(j)M^{2}\Big\}.
\end{eqnarray}
We also have
\begin{equation*}
\mathbb E\left[\left(\sum_{i=1}^{[nt_1]}\vartheta_{ni}^{(1)}\right)\left(\sum_{i=[nt_1]+1}^{n}\vartheta_{ni}^{(2)}\right)\right]
=\sum_{i=1}^{[nt_1]}\sum_{j=[nt_1]+1}^{n}\mathbb E\left(\vartheta_{ni}^{(1)}\vartheta_{nj}^{(2)}\right).
\end{equation*}
From
\begin{equation*}
\sup_{n \ge 1}\sup_{i,j\geq 1}\left\vert \mathbb E\big(\vartheta_{ni}^{(1)}\vartheta_{nj}^{(2)}\big)\right\vert \leq \beta^{\frac{\delta}{2+\delta}}(j-i)MM^{\ast},
\end{equation*}
where $M^{\ast}=\sup_{n \ge 1}\left\{\mathbb E\left[ \left(\vartheta_{n1}^{(2)}\right)^{2+\delta} \right]
\right \}^{\frac{1}{2+\delta}}$, 
it results that
\begin{equation} \label{chp2}
\lim_{n\rightarrow \infty}\frac{1}{n}\sum_{i=1}^{[nt_1]}\sum_{j=[nt_1]+1}^{n}\left\vert
\mathbb E\big(\vartheta_{ni}^{(1)}\vartheta_{nj}^{(2)}\big)\right\vert \leq t_1\sum_{j=1}^{\infty}\beta^{\frac{\delta}{2+\delta}}(j)MM^{\ast}.
\end{equation}
Similarly, we get
\begin{equation} \label{chp3}
\lim_{n\rightarrow
\infty}\frac{1}{n}\mathbb E\left[\bigg(\sum_{i=[nt_1]+1}^{n}\vartheta_{ni}^{(2)}\bigg)^{2}\right]
\leq
(1-t_1)\Big\{\mathbb E\left[\big(\vartheta_{n1}^{(1)}\big)^2 \right]+2\sum_{j=1}^{\infty}\beta^{\frac{\delta}{2+\delta}}(j)(M^{\ast})^{2}\Big\}.
\end{equation}
From (\ref{chp1})-(\ref{chp3}), we deduce (\ref{hp4}).
\\\\
Now, we turn to proving (\ref{hp3}).
For all $i\geq 1$, and for any $\kappa>0$, define
$$
\vartheta_{ni}^{(j),\kappa}=
\left\{
\begin{array}{ll}
\vartheta_{ni}^{(j)} &\mbox{ if } \big|\vartheta_{ni}^{(j)}\big| \leq
\kappa\\
 0 &\mbox{ if } \big|\vartheta_{ni}^{(j)}\big| \geq \kappa, \
j=1,2.
\end{array} \right.
$$
It is immediate that
\begin{equation*}
\sup_{n\geq 1}\sup_{i\geq 1}\big| \vartheta_{ni}^{(j)}\big| \leq \kappa<\infty.
\end{equation*}
It results from the integrability condition in Theorem \ref{th2} that the 
sequences $\{\vartheta_{ni}^{(j)}; \ i\geq 1, \ j=1,2\}$ are uniformly integrable.
\\\\
Whence
\begin{equation*}
\sup_{i\geq 1}\mathbb E\left(\Big| \vartheta_{ni}^{(j)}-\vartheta_{ni}^{(j),\kappa}\Big|^{2+\delta}\right)\longrightarrow 0 \ \text{ as } \ \kappa\rightarrow \infty, \ j=1,2
\end{equation*}
and (\ref{hp3}) is proved.
\\\\
The proof of (\ref{hp5}), that is 
\begin{equation*}
\lim_{n\rightarrow \infty} \mathbb E\left[\left(\frac{1}{\sqrt{n}}\left\{\sum_{i=1}^{[nt_1]}\left[\vartheta_{ni}^{(1),\kappa}-\mathbb E\left(\vartheta_{ni}^{(1),\kappa}\right)\right]
+\sum_{i=[nt_1]+1}^{n}\left[\vartheta_{ni}^{(2),\kappa}-\mathbb E\left(\vartheta_{ni}^{(2),\kappa}\right)\right]\right\}\right)^{2} \right]=c_\kappa, 
\end{equation*}
where $c_\kappa$ is some positive constant, is similar to that of  (\ref{hp4}).
\\\\
It remains to prove (\ref{hp6}).
\medskip

\noindent
For any $i \ge 1$ and any $j=1,2$, denote by $\vartheta_{i}^{(j),\kappa}$ the counterpart of 
$\vartheta_{ni}^{(j),\kappa}$ obtained by substituting the $X_i$'s for the $Y_{ni}^l$'s.
\smallskip

\noindent
We have
\begin{eqnarray*}
c_\kappa & = & t_1 \mathbb E\left\{ \left[\vartheta_{1}^{(1),\kappa}- \mathbb E\big(\vartheta_{1}^{(1),\kappa}\big)\right]^{2} \right\}\\
&&+2t_1\sum_{i=1}^{\infty}
\mathbb E\left\{ \left[\vartheta_{1}^{(1),\kappa}-\mathbb E\big(\vartheta_{1}^{(1),\kappa}\big)\right]\left[\vartheta_{i+1}^{(1),\kappa}-\mathbb E\big(\vartheta_{i+1}^{(1),\kappa}\big)\right]\right\}\\
&  & +t_1\sum_{i=1}^{\infty}\mathbb E\left\{ \left[\vartheta_{1}^{(1),\kappa}-\mathbb E\big(\vartheta_{1}^{(1),\kappa}\big)\right]
\left[\vartheta _{i+1}^{(2),\kappa}-\mathbb E\big(\vartheta_{i+1}^{(2),\kappa}\big)\right]\right\}\\
&  & +(1-t_1) \mathbb E\left\{ \left[\vartheta_{1}^{(2),\kappa}-\mathbb E\big(\vartheta _{1}^{(2),\kappa}\big)\right]^{2} \right\}\\
&  & +2(1-t_1)\sum_{i=1}^{\infty}\mathbb E\left\{ \left[\vartheta _{1}^{(2),\kappa}-\mathbb E\big(\vartheta _{1}^{(2),\kappa}\big)\right]\left[\vartheta_{i+1}^{(1),\kappa}-\mathbb E\big(\vartheta _{i+1}^{(2),\kappa}\big)\right]\right\}.
\end{eqnarray*}
By the Lebesgue dominated convergence theorem, one obtains
\begin{equation*}
\mathbb E\left\{ \left[\vartheta _{1}^{(1),\kappa}-\mathbb E\left(\vartheta _{1}^{(1),\kappa}\right)\right]^{2} \right \}\longrightarrow
\mathbb E\left[\left(\vartheta_{1}^{(1)}\right)^{2} \right] \mbox{ \ as } \kappa \rightarrow \infty,
\end{equation*}
\begin{equation*}
\mathbb E\left\{ \left[\vartheta_{1}^{(1),\kappa}-\mathbb E\big(\vartheta_{1}^{(1),\kappa} \right] \left[(\vartheta_{i+1}^{(1),\kappa}-\mathbb E\big(\vartheta_{i+1}^{(1),\kappa}\big)\right]\right \} \longrightarrow \mathbb E \left(\vartheta_{1}^{(1)}\vartheta_{i+1}^{(1)}\right) \mbox{ \ as }\kappa\rightarrow \infty,
\end{equation*}
\begin{equation*}
\mathbb E\left\{ \left[\vartheta_{1}^{(1),\kappa}-\mathbb E\big(\vartheta_{1}^{(1),\kappa}\big)\right]\left[\vartheta_{i+1}^{(2),\kappa}-\mathbb E\big(\vartheta_{i+1}^{(2),\kappa}\big)\right] \right\}\longrightarrow \mathbb E\left(\vartheta_{1}^{(2)}\vartheta_{i+1}^{(2)}\right) \mbox{ \ as } \kappa\rightarrow \infty,
\end{equation*}
\begin{equation*}
\mathbb E\left\{ \left[\vartheta_{1}^{(2),\kappa}-\mathbb E\big(\vartheta_{1}^{(2),\kappa}\big)\right]^{2}\right\} \longrightarrow \mathbb E\left[\left(\vartheta_{1}^{(2)}\right)^{2} \right] \mbox{ \ as } \kappa\rightarrow \infty
\end{equation*}
and
\begin{equation*}
\mathbb E\left\{ \left[\vartheta_{1}^{(2),\kappa}-\mathbb E\big(\vartheta_{1}^{(2),\kappa}\big)\right]\left[\vartheta_{i+1}^{(2),\kappa}
-\mathbb E\big(\vartheta_{i+1}^{(2),\kappa}\big)\right] \right \} \longrightarrow \mathbb E\left( \vartheta_{1}^{(2)} \vartheta_{i+1}^{(2)} \right) \mbox{ \ as } \kappa\rightarrow \infty.
\end{equation*}
Therefore
\begin{equation*}
\lim_{\kappa \rightarrow \infty }c_\kappa=c
\end{equation*}
and (\ref{hp6}) is proved.
Whence, 
the finite dimensional convergence is established.
\bigskip

\noindent
For proving the tightness define
\begin{equation*}
\mathcal{Q}_n(t^l)=\sigma_{11}^{-1/2}\frac{1}{\sqrt n}\sum_{i=1}^{[nt^l]}h_{F_l^{(n)},1}(Y^l_{ni}).
\end{equation*}
If $t_1^l\leq t^l\leq t_2^l$, from the integral conditions and condition (\ref{bet}), there exists a constant $C$ such that
\begin{eqnarray*}
\mathbb{E}\left(\left|\mathcal{Q}_n(t^l)-\mathcal{Q}_n(t_1^l)\right|^2
\left|\mathcal{Q}_n(t_2^l)-\mathcal{Q}_n(t^l)\right|^2\right) & \leq &
C\frac{1}{n^2}\left([nt^l]-[nt_1^l]\right)\left([nt_2^l]-[nt^l]\right)\\
& \leq & C\frac{1}{n^2}\left([nt_1^l]-[nt^l]\right)\left([nt^l]-[nt_2^l]\right)\\
& \leq & C\frac{1}{n^2}\left([nt_2^l]-[nt_1^l]\right)^2\\
& \leq & C\left(t_2^l-t_1^l\right)^2.
\end{eqnarray*}
If $t_2^l-t_1^l\geq 1/n$ the last inequality follows and if $t_2^l-t_1^l< 1/n$,
then either $t_1^l$ and $t^l$ lie in the same subinterval $[(i-1)/n,i/n]$ or else $t^l$ and $t_2^l$ do.
In either of these cases the left hand of last inequality vanishes. From
Theorem 13.5 of \cite{billingsley}, the process $\mathcal{Q}_n$ is tight.
\medskip

\noindent
With similar methods, we prove that any linear combinations of the components converges to a Gaussian random variable. 
Therefore Proposition \ref{pro3} is proved by Lemma \ref{phil}.
\bigskip

\subsection{Proof of Theorem \ref{th1}}
Using the Hoeffding decomposition, for any $(t_1, \ldots,t_k) \in \Theta_k$, we can write $\mathcal Z_n(t_1,\ldots,t_k)$ as
\begin{eqnarray}
&&\mathcal Z_{n}(t_1,\ldots,t_k) = n^{-3/2}\sum_{l=1}^{k}\sum_{i=[nt_{l-1}]+1}^{[nt_l]}\sum_{j=[nt_l]+1}^{[nt_{l+1}]}
\left[h_{F_1,1}(X_{i})+h_{F_1,2}(X_{j})+g_{F_1,F_1}(X_{i},X_{j})\right]
\nonumber\\
 & = & n^{-3/2}\sum_{l=1}^{k}\left[([nt_{l+1}]-[nt_l])\sum_{i=[nt_{l-1}]+1}^{[nt_l]}h_{F_1,1}(X_{i})+([nt_l]-[nt_{l-1}])
 \sum_{j=[nt_l]+1}^{[nt_{l+1}]}h_{F_1,2}(X_{j})\right]\nonumber\\
 & & +n^{-3/2}\sum_{l=1}^{k}\sum_{i=[nt_{l-1}]+1}^{[nt_l]}\sum_{j=[nt_l]+1}^{[nt_{l+1}]}g_{F_1,F_1}(X_{i},X_{j}), \label{z1}
\end{eqnarray}
where by convention $[nt_0]=1$ and $[nt_{k+1}]=n$.
\medskip

From Proposition \ref{pro3}, as $n$ tends to infinity, in probability, we have,
\begin{equation*}
n^{-3/2}\sup_{(t_1, \ldots,t_k) \in \Theta_k}
\left|\sum_{l=1}^{k}\sum_{i=[nt_{l-1}]+1}^{[nt_l]}\sum_{j=[nt_l]+1}^{[nt_{l+1}]}
g_{F_1,F_1}(X_{i},X_{j})\right|\longrightarrow 0.
\end{equation*}

\noindent
Thus, by Slutsky's lemma, it suffices to show that the sequence of the following processes indexed by $(t_1,\ldots,t_k) \in \Theta_k$
\begin{equation*}
n^{-3/2}\sum_{l=1}^{k}\left[([nt_{l+1}]-[nt_l])\sum_{i=[nt_{l-1}]+1}^{[nt_l]}h_{F_1,1}(X_{i})+([nt_l]-[nt_{l-1}])
 \sum_{j=[nt_l]+1}^{[nt_{l+1}]}h_{F_1,2}(X_{j})\right]
\end{equation*}
converges weakly to the desired limit process.
\medskip

\noindent
Now, observing that the following mapping 
\begin{equation*}
\begin{pmatrix}
   x_{1}(t_1) \\
   x_{2}(t_1) \\
   x_{3}(t_2) \\
   x_{4}(t_2) \\
   \vdots \\
   x_{l}(t_l) \\
   x_{l+1}(t_l) \\
   \vdots \\
   x_{k}(t_k) \\
   x_{k+1}(t_k) \\
\end{pmatrix}%
\longmapsto \sum_{l=1}^{k} \{(t_{l+1}-t_l)[x_{l}(t_{l})-x_{l}(t_{l-1})] 
+(t_l-t_{l-1})[x_{l+1}(t_{l+1})-x_{l}(t_l)]\}
\end{equation*}
 is continuous from $(D[0,1])^{2k}$ to $D[0,1]$, one has by the continuous mapping theorem and Proposition \ref{pro2}, the following weak convergence
\begin{eqnarray*}
&&n^{-3/2}\displaystyle\sum_{l=1}^{k}([nt_{l+1}]-[nt_l])\sum_{i=[nt_{l-1}]+1}^{[nt_l]}h_{F_1,1}(X_{i})\\&&+n^{-3/2}\sum_{l=1}^{k}([nt_l]-[nt_{l-1}])
\displaystyle\sum_{j=[nt_l]+1}^{[nt_{l+1}]}h_{F_1,2}(X_{j})
% \hspace{5cm}&\\
%&\hspace{8cm}
\\&& \longrightarrow \mathcal Z(t_1,\ldots,t_k)
\end{eqnarray*}
and Theorem \ref{th1} is proved. \hspace{0.5cm} $\Box$

\subsection{Proof of Theorem \ref{th2}}
For the proof of Theorem \ref{th2}, we observe that 
under its assumptions,
% letting  
%Under $\mathcal{H}_1^k$, for $1\leq l\leq k+1$, consider the following
for the  random variables $V_{n,\dot F_l}$ defined by 
$$V_{n,\dot F_l}=\sum_{i=[nt_{l-1}]+1}^{[n t_l]}
\sum_{j=[nt_l]+1}^{[nt_{l+1}]}h(X_{i},X_{j}),$$
one can check easily that 
\begin{eqnarray} \label{eq5}
&& V_{n,\dot F_l}-([nt_{l}]-[nt_{l-1}])([nt_{l+1}]-[nt_{l}])\theta_h(\dot F_{l-1},\dot F_{l-1})  \nonumber\\
& =& ([nt_{l}]-[nt_{l-1}])([nt_{l+1}]-[nt_{l}])\theta_h(\dot F_{l-1},\dot F_{l})-\theta_h(\dot F_{l-1},\dot F_{l-1}) \nonumber\\
& +&([nt_{l+1}]-[nt_{l}])\sum_{i=[nt_{l-1}]+1}^{[nt_{l}]}h_{\dot F_{l},1}(X_{i})
+([nt_{l}]-[nt_{l-1}])\sum_{i=[nt_{l}]+1}^{[nt_{l+1}]}h_{\dot F_{l},2}(X_{i})  \nonumber\\
&+& \sum_{i=[nt_{l-1}]+1}^{[nt_{l}]}\sum_{j=[nt_{l}]+1}^{[nt_{l+1}]}g_{\dot F_{l-1},\dot F_l}(X_{i},X_{j}).
\end{eqnarray}

Thus, it results from (\ref{eq5}), that for any $(t_1, \ldots,t_k) \in \Theta_k$, the following equalities hold 
\begin{eqnarray*}
&&\mathcal Z_{n}(t_1,\ldots,t_k) = n^{-3/2}\sum_{l=1}^kV_{n,F^{(n)}_l}\\
& & -n^{-3/2}\sum_{l=1}^k([nt_{l}]-[nt_{l-1}])([nt_{l+1}]-[nt_{l}])\left[\theta_h(F^{(n)}_{l-1},F^{(n)}_l)-\theta_h(F^{(n)}_{l-1},F^{(n)}_{l-1}) \right]\\
& = & n^{-3/2}\sum_{l=1}^k\bigg\{([nt_{l+1}]-[nt_{l}])\sum_{i=[nt_{l-1}]+1}^{[nt_{l}]}h_{F^{(n)}_{l},1}(Y^{l-1}_{ni})
\\ & & 
+([nt_{l}]-[nt_{l-1}])\sum_{j=[nt_{l}]+1}^{[nt_{l+1}]}h_{F^{(n)}_{l},2}(Y^l_{nj})\\
& & +\sum_{i=[nt_{l-1}]+1}^{[nt_{l}]}\sum_{j=[nt_{l}]+1}^{[nt_{l+1}]}g_{F^{(n)}_{l-1},F^{(n)}_l}(Y^{l-1}_{ni},Y^l_{nj})\bigg\}\\
& & -n^{-3/2}\sum_{l=1}^k([nt_{l}]-[nt_{l-1}])([nt_{l+1}]-[nt_{l}])\left[\theta_h(F^{(n)}_{l-1},F^{(n)}_l)-\theta_h(F^{(n)}_{l-1},F^{(n)}_{l-1}) \right].
\end{eqnarray*}
From Proposition \ref{pro2}, we deduce that as $n$ tends to infinity, in probability, 
\begin{equation*}
n^{-3/2}\sup_{(t_1, \ldots,t_k) \in \Theta_k}
\left|\sum_{l=1}^k\sum_{i=[nt_{l-1}]+1}^{[nt_l]}\sum_{j=[nt_l]+1}^{[nt_{l+1}]}g_{F^{(n)}_{l-1},F^{(n)}_l}(Y^{l-1}_{ni},Y^l_{nj})\right|
\longrightarrow 0.
\end{equation*}

\noindent
By Proposition \ref{pro2}, it is easy to see that, for any $l=1, \ldots,k$, the sequence of processes 
\begin{equation*}
n^{-1/2}\sum_{i=1}^{[nt_l]}h_{F^{(n)}_{l},1}(Y^{l-1}_{ni})
\end{equation*}
converges weakly to the Brownian process $\{W_1(t_l)\}_{_{0\leq t_l \leq 1}}$ and that the sequence of processes 
\begin{equation*}
n^{-1/2}\sum_{j=1}^{[nt_{l}]}h_{F_l^{(n)},2}(Y_{nj}^l)
\end{equation*}
converges weakly to the Brownian process $\{W_2(t_l)\}_{_{0\leq t_l \leq 1}}$.

\noindent
As in the proof of Theorem \ref{th1}, there is a continuous mapping from $(D[0,1])^{2k}$ to $D[0,1]$ such that, as $n$ tends to infinity, the following weak convergence holds 
\begin{eqnarray*}
&n^{-3/2}\displaystyle\sum_{l=1}^k\bigg\{([nt_{l+1}]-[nt_{l}])\sum_{i=[nt_{l-1}]+1}^{[nt_{l}]}h_{F^{(n)}_{l},1}(Y^{l-1}_{ni})
+([nt_{l}]-[nt_{l-1}])\sum_{j=[nt_{l}]+1}^{[nt_{l+1}]}h_{F^{(n)}_{l},2}(Y^l_{nj})\bigg\}&\\
&\hspace{7cm} \longrightarrow \mathcal Z(t_1,\ldots,t_k).&
\end{eqnarray*}
Also, under $H_{1,n}^k$, as $n$ tends to infinity,  

\begin{eqnarray*}
&n^{-3/2}\displaystyle\sum_{l=1}^k([nt_{l}]-[nt_{l-1}])([nt_{l+1}]-[nt_{l}]) \left[\theta_h(F^{(n)}_{l-1},F^{(n)}_l)-\theta_h(F^{(n)}_{l-1},F^{(n)}_{l-1}) \right] \hspace{2cm} &\\
&\hspace{7cm} \longrightarrow \displaystyle \sum_{l=1}^k(t_{l+1} -t_l)(t_l-t_{l+1})A_l&
\end{eqnarray*}
and Theorem \ref{th2} is proved. \hspace{0.5cm} $\Box$

\subsection{Proof of Theorem \ref{th3}}
Let $[nt_{0,l-1}]\le[nt_l]\le[nt_{0,l}], 1\le l \le k$, then
\begin{eqnarray*}
\mathcal Z_{n}^*(t_1,\ldots,t_k) & = & n^{-3/2}\sum_{l=1}^k\sum_{[nt_{0,l-1}]\le i < j \le [nt_{0,l}]}h(X_i,X_j)
\\&&+ n^{-3/2}\sum_{l=1}^k\sum_{i=[nt_{l-1}]+1}^{[nt_{0,l}]}\sum_{j=[nt_{0,l}]+1}^{[nt_{0,l+1}]}h(X_i,X_j)\\
&  & -\bigg\{n^{-3/2}\sum_{l=1}^k\sum_{[nt_{0,l-1}] \le i<j \le [nt_l]}
h(X_i,X_j)\\
&&
+n^{-3/2}\sum_{l=1}^k\sum_{[nt_l]+1\leq i<j\leq [nt_{0,l}]}h(X_i,X_j)\\
&  & +n^{-3/2}\sum_{l=1}^k\sum_{[nt_l]+1\leq i\leq [nt_{0,l}]}\sum_{[nt_{0,l}]+1\leq j\leq [nt_{0,l+1}]}h(X_i,X_j)\bigg\} \\
& = & R_n^{(1)} + R_n^{(2)} - \left\{R_n^{(3)} + R_n^{(4)} + R_n^{(5)}\right\}.
\end{eqnarray*}
Similar to the proof of Theorem 3 (see \cite{nga}).\\
First we prove that as $n$ tends to infinity, almost surely, 
\begin{equation*}
n^{-1/2}R_n^{(1)} \longrightarrow \sum_{l=1}^k{t_{0,l}}^2\theta_h(\dot F_l, \dot F_l)/2.
\end{equation*}
Similarly, we prove that as $n$ tends to infinity, almost surely, 
\begin{equation*}
n^{-1/2}R_n^{(3)} \longrightarrow  \sum_{l=1}^kt_l^2\theta_h(\dot F_l, \dot F_l)/2,
\end{equation*}
\begin{equation*}
n^{-1/2}R_n^{(4)}\stackrel{ \mathcal{D} }{=}\sum_{l=1}^k\sum_{1\leq i < j \leq[nt_{0,l}] -[(n+1)t_l]}h(X_i,X_j) \longrightarrow \sum_{l=1}^k(t_l-t_{0,l})^2\theta_h(\dot F_l,\dot F_l)/2.
\end{equation*}
and, we establish that, as $n$ tends to infinity, almost surely, 
\begin{equation*}
n^{-1/2}R_n^{(2)} \longrightarrow \sum_{l=1}^kt_{0,l}(t_{0,l+1}-t_{0,l})\theta_h(\dot F_l,\dot F_{l+1}).
\end{equation*}
Similarly, we prove that, as $n$ tends to infinity, in probability, 
\begin{equation*}
n^{-1/2}R_n^{(5)} \longrightarrow \sum_{l=1}^k(t_{0,l}-t_l)(t_{0,l+1}-t_{0,l})\theta_h(\dot F_l,\dot F_{l+1}).
\end{equation*}
The first part of (\ref{z_1*}) clearly follows from these results. The proof of its second part can be handled in similar lines.  \hspace{0.5cm} $\Box$

%\bibliographystyle{apalike}
%\bibliography{RuptureMultiple}
%\bibliography{Bibliography-jasa_es1}
%\bibliographystyle{apalike}
%\end{document}
%

\end{document}